\documentclass[10pt]{article}
\usepackage{amsmath,amssymb,amsthm,float,bbm}
\usepackage{xcolor}
\usepackage{tikz-cd}
\usepackage[margin=1in]{geometry}
\usepackage[hypertexnames=false]{hyperref}
\theoremstyle{plain}
\newtheorem{theorem}{Theorem}[section]
\newtheorem{proposition}[theorem]{Proposition}
\newtheorem{lemma}[theorem]{Lemma}
\newtheorem{corollary}[theorem]{Corollary}
\theoremstyle{definition}
\newtheorem{definition}[theorem]{Definition}
\newtheorem{note}[theorem]{Note}
\newtheorem{example}[theorem]{Example}
\theoremstyle{plain}
\newtheorem{innercustomthm}{Theorem}
\newenvironment{customthm}[1]{\renewcommand\theinnercustomthm{#1}\innercustomthm}{\endinnercustomthm}

\newcommand{\Z}{\mathcal{Z}}
\newcommand{\Jac}[2]{\left(\tfrac{#1}{#2}\right)}

\title{Stretched near-groups}
\author{Andrew Schopieray, Yilong Wang, Zhiqiang Yu}
\date{\today}

\begin{document}
\maketitle

\begin{abstract}
\noindent We study stretched near-group fusion categories over $\mathbb{C}$ of type $G+k|G|$ for abelian groups $G$ and positive integers $k$, describing the structure of the double including fusion rules, modular data, and Frobenius--Schur indicators under a variety of assumptions. For cyclic $G$, when enough invertible objects exist in the double, we describe all associated de-equivariantizations which connect the existence of stretched near-group fusion categories to the existence of other stretched quadratic categories. All of our results apply to the case $k=1$. In particular, the known existence of $C_{m^2}+m^2$ near-group categories for all positive integers $m$ implies the existence of Haagerup--Izumi fusion categories for all cyclic groups. As a further application, we provide a general conductor bound severely restricting the existence of stretched near-group categories, and for infinitely many abelian groups $G$ of square order, among them $C_4\times C_2^{2j}$ with $j\geq2$, we show that a near-group fusion category $G+k|G|$ exists if and only if $k=0$.
\end{abstract}
\tableofcontents

\section{Introduction}

\par In this paper, we study near-group fusion categories of type $G+k|G|$ with $G$ a finite abelian group and $k\in\mathbb{Z}_{\geq1}$, i.e.\ fusion categories with a single isomorphism class $x$ of noninvertible simple objects such that $x\otimes x$ contains $x$ with multiplicity $k|G|$. We will refer to the underlying fusion rings as \emph{stretched} near-group fusion rings; for $k=1$ these are the near-groups of type $G+|G|$, which have been studied extensively but for the most part under strong operator algebraic assumptions (see for example \cite{MR2837122,MR3167494,MR3635673,MR5015436}).  Renewed interest in these topics came from a proof of D.\ Radchenko and C.\ Wheeler that $C_n\times C_m+nm$ near-group fusion categories exist for all $n,m\geq1$ \cite{radchenko2026realquadraticfieldsfinite} which is closely related to the existence proof of Haagerup--Izumi categories for odd-order cyclic groups due to T.\ Gannon, H.\ Yadav, and the first author \cite{gannon2026haagerupizumifusioncategories} and independently T.-C.\ Huang \cite{huang2026cyclichaagerupizumifusioncategories}.  These results all tie the existence of such fusion categories to the algebraicity of values of special functions such as quantum dilogarithms and double-sines. With these deep connections between number theory and quadratic fusion categories, understanding the structure of stretched near-group categories has become a vital but seemingly monumental task.

\par Our main structural results surround the doubles $\mathcal{Z}(\mathcal{N}_G)$ of stretched near-group fusion categories $\mathcal{N}_G$. The pointed part $(\mathcal{N}_G)_\mathrm{pt}\simeq\mathrm{Vec}_G$ corresponds to a connected \'etale algebra $A_G$ in the double whose category of local modules is the untwisted double $\mathcal{Z}(\mathrm{Vec}_G)$, and the branching rules control the inductions $I(g)$ of the invertible objects. Section \ref{sec:prelim} collects conventions, a description of inductions as local modules over the connected \'etale subalgebras of the canonical Lagrangian algebra $I(\mathbbm{1})$ attached to fusion subcategories (Proposition \ref{correspond}), a bound on the twists of inductions in arbitrary spherical fusion categories (Proposition \ref{twist}), and a multiplicativity property of Frobenius--Schur indicators for products with central objects (Lemma \ref{lem:nucentral}). In Section \ref{sec:neargroup} we show that each $I(g)$ is multiplicity free: it is a sum of $|G|$ objects, all simple except one, which splits into two simple objects $Y_g$ and $Y_g'$ (Proposition \ref{prop:local}). From the twists of these summands we construct, for every $k$, the expected (as with $k=1$) quadratic form $q$ on $G$ with associated bicharacter $b$ (Proposition \ref{prop:tq}). When $G$ is cyclic, this form can be determined completely (Proposition \ref{cor:cyclic}).

\par Section \ref{sec:neargroup} also treats the higher Frobenius--Schur indicators of $x$, for every $k$ and without further hypotheses. For every $r$ coprime to the Frobenius--Schur exponent $N$, $\nu_r(x)=\frac{k}{2}(1-\Jac{\mathfrak{d}}{r})$, where $\mathfrak{d}$ is the discriminant of $\mathbb{Q}(\sqrt{|G|(k^2|G|+4)})$ (Theorem \ref{thm:nugeneral}). Moreover, the primes dividing $N$ are exactly those dividing $|G|(k^2|G|+4)$, and $N$ is even whenever $k\geq2$. In particular, if $|G|$ is odd, then $G+k|G|$ is categorifiable only if $k=1$ or $k$ is even (Corollary \ref{cor:fsexp}). The indicators constrain the twists of the simple objects of the double that do not lie over invertible objects, through a family of identities (Lemma \ref{lem:moments}).

\par For $k=1$, and for every $k$ under the hypothesis \eqref{eq:star} that $\mathcal{O}(\mathcal{Z}(\mathcal{N}_G)_\mathrm{pt})\cong G$, the radical of $b$ is the $2$-torsion subgroup $G[2]$ (Theorem \ref{thm:radical}). In particular, $\mathcal{Z}(\mathcal{N}_G)\simeq\mathcal{Z}(\mathcal{N}_G)_\mathrm{pt}\boxtimes\mathcal{Z}(\mathcal{N}_G)_\mathrm{ad}$ if and only if $|G|$ is odd (Corollaries \ref{cor:oddfactor} and \ref{cor:factiff}). Motivated by this, we call $\mathcal{N}_G$ \emph{$G$-factorizable} if its double factors in this way (Definition \ref{def:factorizable}). In Section \ref{sec:factorizable} we determine, for $G$-factorizable near-groups of type $G+k|G|$, the simple objects of $\mathcal{Z}(\mathcal{N}_G)_\mathrm{ad}$ (Proposition \ref{factorprop}), most of its fusion rules (Propositions \ref{fus}--\ref{prop:yeprime}), and most of its modular data (Section \ref{subsec:modulardata}). The remaining simple objects $Z_j$, with forgetful images $m_jx$, satisfy $m_j^2\mid k^2|G|+4$; in particular, when $k\geq3$ none of them lies over $kx$ (Corollary \ref{cor:codeg}), whereas for $k=1$ all of them lie over $x$ whenever no square greater than $1$ properly divides $|G|+4$, for instance when $|G|+4$ is squarefree (Corollary \ref{cor:codeg} and Note \ref{rem:codeg}). The part of the $T$-matrix that remains unknown is constrained by the twist identities, which in this case reduce to one identity for each $r\geq1$ involving the Gauss sums of $\mathcal{C}(G,q)$ (Theorem \ref{thm:master}); for $r$ coprime to $N$ these identities are Galois conjugate to the one for $r=1$ (Note \ref{rem:nunGalois}).

\par In Section \ref{sec:deeq} we consider $G=C_n$ under hypothesis \eqref{eq:star}: that $\mathcal{O}(\mathcal{Z}(\mathcal{N}_G)_\mathrm{pt})\cong G$. Writing $n=m^2s$ with $s$ squarefree, the maximal Tannakian subcategory of $\mathcal{Z}(\mathcal{N}_G)_\mathrm{pt}$ is $\mathrm{Rep}(C_m)$, independently of $q$ and $k$ (Proposition \ref{prop:maxtann}). For odd $n>1$ it is Lagrangian exactly when $n$ is a perfect square, in which case $\mathcal{Z}(\mathcal{N}_G)_\mathrm{pt}\simeq\mathcal{Z}(\mathrm{Vec}_{C_m}^\omega)$ with $\omega$ of order $m$ (Corollary \ref{cor:lag}). We then describe the de-equivariantizations by these Tannakian subcategories (Proposition \ref{prop:deeq}). When $n=m^2$ is a perfect square, the result categorifies a \emph{stretched Haagerup--Izumi fusion ring} $H_k^z(C_m)$ for some $z\in C_m$ with $z^2=e$: the fusion ring with basis $\{g,\,g\rho:g\in C_m\}$ and $\rho g=g^{-1}\rho$, $\rho^2=z+k\sum_{g\in C_m}g\rho$. We write $H_k(C_m):=H_k^e(C_m)$, in which every $g\rho$ is self-dual; if $z\neq e$, then $z$ is the element of order $2$ and no $g\rho$ is self-dual.

\begin{customthm}{\ref{thm:haagerup}}
Let $m>1$, and let $\mathcal{N}_G$ be a near-group fusion category of type $C_{m^2}+km^2$ satisfying \eqref{eq:star}, for instance $k=1$. Then the maximal Tannakian subcategory of $\mathcal{Z}(\mathcal{N}_G)_\mathrm{pt}$ is $\mathrm{Rep}(C_m)$, and the de-equivariantization $\mathcal{D}=(\mathcal{N}_G)_{C_m}$ satisfies $K(\mathcal{D})\cong H_k(C_m)$ if $m$ is odd and $K(\mathcal{D})\cong H_k^z(C_m)$ if $m$ is even, where $z\in C_m$ is the element of order $2$.
\end{customthm}

\par Combined with the existence of $C_{m^2}+m^2$ near-group fusion categories for every $m\geq1$ \cite{radchenko2026realquadraticfieldsfinite}, Theorem \ref{thm:haagerup} provides, for every $m>1$, a fusion category categorifying $H_1(C_m)$ if $m$ is odd and $H_1^z(C_m)$ if $m$ is even. In each case the invertible objects form a subcategory equivalent to $\mathrm{Vec}_{C_m}^\omega$ with $\omega$ of order $m$ (Proposition \ref{prop:Dpt}). For $m$ odd these are therefore \emph{twisted} Haagerup--Izumi fusion categories, as opposed to the infinite families of odd Haagerup--Izumi categories constructed in \cite{gannon2026haagerupizumifusioncategories}
 and \cite{huang2026cyclichaagerupizumifusioncategories} which have pointed subcategory $\mathrm{Vec}_{C_m}$. For $m$ even they form an infinite family of Haagerup--Izumi fusion categories, in the sense of Definition \ref{def:haagerup}, for cyclic groups of \emph{even} order, none of whose noninvertible simple objects is self-dual. This goes beyond the existence results of \cite{gannon2026haagerupizumifusioncategories} and \cite{huang2026cyclichaagerupizumifusioncategories}, which concern Haagerup--Izumi fusion rings for cyclic groups of odd order. The proof for even $m$ uses no explicit construction: the element of order $2$ in $\mathcal{D}$ has a central lift whose twist is a primitive fourth root of unity, so its Frobenius--Schur indicator is $-1$, and Frobenius--Schur indicators are multiplicative against central objects (Lemma \ref{lem:nucentral}). We note that these results would follow from the work of Izumi \cite[Section 12]{MR3635673} for a unitary version of the Radchenko-Wheeler near-group constructions \cite[Section 5]{radchenko2026realquadraticfieldsfinite} but no such solutions are known at this time.

\par Finally, Section \ref{sec:application} treats arbitrary categorifications of $G+k|G|$, without \eqref{eq:star} or $G$-factorizability. We prove a conductor bound coming from the dimension block of a spherical fusion category in the sense of \cite{schopieray2026dimensionblocksphericalfusion}: if $G+k|G|$ is categorifiable, then $k^2\mathfrak{f}\leq|G|(k^2|G|+4)$, where $\mathfrak{f}$ is the conductor of $\mathbb{Q}(\sqrt{|G|(k^2|G|+4)})$ (Proposition \ref{prop:conductor}). This bound is sharp for $C_3+6$. Section \ref{subsec:2groups} then combines Corollary \ref{cor:fsexp} with the twist identities of Lemma \ref{lem:moments} and averages over Galois groups (Lemma \ref{lem:average}). Write $R(G,\ell)$ for the near-group fusion ring $x^2=\sum_{g\in G}g+\ell x$. The levels $\ell$ for which $R(G,\ell)$ is categorifiable were, to our knowledge, previously known only for $G=C_1$, $C_3$ and the elementary abelian $2$-groups (including $C_2$). We determine them for infinitely many abelian groups of square order, using Larson-type bounds on sums of roots of unity (Lemma \ref{lem:average}(2)). Write $2G:=\{g^2:g\in G\}$.

\begin{customthm}{\ref{thm:2groups}}
For every $a\in\mathbb{Z}_{\geq1}$ there is a positive integer $j_0(a)$ with the following property. Let $G$ be an abelian group of order $4^jw^2$ with $w$ odd, $|2G|=a$ and $j\geq j_0(a)$. Then $G+k|G|$ is not categorifiable for any $k\geq1$, and $R(G,\ell)$ is categorifiable if and only if $\ell=0$. For example, one can take $j_0(2)=3$, that is, the conclusion holds for $G\cong C_4\times C_2^{2j-2}$ with $j\geq3$.
\end{customthm}


\section{Preliminaries}\label{sec:prelim}


\par We will only discuss fusion categories over $\mathbb{C}$, in the sense of \cite{MR2183279}: semisimple tensor categories $\mathcal{C}$ with finitely many isomorphism classes of simple objects, which we denote $\mathcal{O}(\mathcal{C})$, and a simple unit $\mathbbm{1}$. Our main tool is the adjoint pair of functors $I:\mathcal{C}\to\mathcal{Z}(\mathcal{C})$ and $F:\mathcal{Z}(\mathcal{C})\to\mathcal{C}$ between $\mathcal{C}$ and its Drinfeld center, or \emph{double}, $\mathcal{Z}(\mathcal{C})$: a braided fusion category whose objects are objects of $\mathcal{C}$ equipped with half-braidings. When $\mathcal{C}$ is spherical, $\mathcal{Z}(\mathcal{C})$ is a modular fusion category. We write $\oplus$ and $+$ interchangeably when we study the fusion rings.


\subsection{Conventions}\label{subsec:conventions}

\par For a spherical fusion category $\mathcal{C}$ we write $\mathrm{Tr}$ for the (spherical) trace, $\dim(X)=\mathrm{Tr}(\mathrm{id}_X)$, and $\dim(\mathcal{C})=\sum_{X\in\mathcal{O}(\mathcal{C})}\dim(X)^2$. The twist of a simple object $Y$ of the modular fusion category $\mathcal{Z}(\mathcal{C})$ is denoted $\theta_Y$, and the $T$-matrix is the diagonal matrix of twists; its order is the Frobenius--Schur exponent $\mathrm{FSexp}(\mathcal{C})$ \cite[Theorem 5.5]{MR2313527}. For a modular fusion category we write $S_{X,Y}$ for the (unnormalized) $S$-matrix, $N_{X,Y}^W=[W,X\otimes Y]$ for the fusion coefficients, where in general $[X,Y]:=\dim\mathrm{Hom}(X,Y)$ (for $Y$ simple, the multiplicity of $Y$ in $X$; this is symmetric in $X$ and $Y$), and we use the balancing equation in the form
\begin{equation}\label{eq:balancing}
    \theta_X\theta_YS_{X,Y}=\sum_{W}N_{X^\ast,Y}^W\,\theta_W\dim(W).
\end{equation}
For $r\in\mathbb{Z}_{\geq1}$ and $X\in\mathcal{C}$, $\nu_r(X)$ denotes the $r$-th Frobenius--Schur indicator of $X$ \cite{MR2313527}; in particular $\nu_1(X)=\delta_{X,\mathbbm{1}}$ for simple $X$. The indicators may be computed from the double via
\begin{equation}\label{eq:NS}
    \nu_r(X)\dim(\mathcal{C})=\mathrm{Tr}\bigl(\theta^r_{I(X)}\bigr)
\end{equation}
by \cite[Theorem 4.1]{MR2313527}. Every simple summand of $I(\mathbbm{1})$ has trivial twist \cite[Remark 4.6]{MR2313527}, \cite[Theorem 2.5]{MR3427429}. We also set $\nu_r(\mathcal{C}):=\sum_{X\in\mathcal{O}(\mathcal{C})}\dim(X)\nu_r(X)$.
Finally, for a (pre-)modular fusion category $\mathcal{B}$, $\tau_r(\mathcal{B})=\sum_{X\in\mathcal{O}(\mathcal{B})}\theta_X^r\dim(X)^2$ denotes its $r$-th Gauss sum \cite{MR3997136}.

\subsection{The correspondence theorem}\label{app:correspond}

Recall \cite[Definition 3.12]{MR3039775} that for a connected \'etale algebra $A$ in a
braided fusion category $\mathcal{A}$, a right $A$-module $M \in \mathcal{A}_A$ is
\emph{dyslectic}, or \emph{local}, if $\rho \circ c_{A,M}\circ c_{M,A} = \rho$,
where $\rho \colon M \otimes A \to M$ is the action; the local $A$-modules form a
fusion subcategory $\mathcal{A}^0_A \subseteq \mathcal{A}_A$, which inherits a braiding from $\mathcal{A}$. We recall the notation of the proof of \cite[Theorem 4.10]{MR3039775}. For a fusion category $\mathcal{C}$ and a fusion subcategory $\mathcal{B}\subset\mathcal{C}$, the
\emph{relative center} $\Z_\mathcal{B}(\mathcal{C})$ is the tensor category of pairs
$(X,\gamma_X)$, where $X\in\mathcal{C}$ and
$\gamma_X:Y\otimes X\stackrel{\sim}{\to}X \otimes Y$ for $Y\in\mathcal{B}$, is
a natural family of isomorphisms satisfying the same compatibility condition as
in the definition of $\Z(\mathcal{C})$.

\par The forgetful functor $F:\Z(\mathcal{C})\to\mathcal{C}$ factors as $\Z(\mathcal{C})\stackrel{F_\mathcal{B}}{\longrightarrow}\Z_\mathcal{B}(\mathcal{C})\stackrel{\widetilde{F}_\mathcal{B}}{\longrightarrow}\mathcal{C}$ where $F_\mathcal{B}$ and $\widetilde{F}_\mathcal{B}$ are the obvious forgetful functors (forgetting some of the half-braidings and forgetting the rest); let $I_\mathcal{B}$
and $\widetilde{I}_\mathcal{B}$ denote their right adjoints so that the adjoint of $F$ is the induction to the center $I=I_\mathcal{B}\circ\widetilde{I}_\mathcal{B}$. The connected \'etale algebra associated with $\mathcal{B}$ in \cite[Theorem 4.10]{MR3039775} is $A:=I_\mathcal{B}(\mathbbm{1})\subset I(\mathbbm{1})$.

\begin{proposition}\label{correspond}
Let $\mathcal{C}$ be a fusion category, $\mathcal{B}\subseteq\mathcal{C}$ a fusion subcategory, and
$A\subset I(\mathbbm{1})$ the connected \'etale algebra
corresponding to $\mathcal{B}$.
\begin{enumerate}
\item For every $X\in\mathcal{O}(\mathcal{B})$, the object $I(X)$ admits the structure of a local $A$-module, i.e.\ it is the underlying object of an object of $\mathcal{Z}(\mathcal{C})_A^0$. More precisely, under the braided equivalence $\mathcal{Z}(\mathcal{C})_A^0\simeq\mathcal{Z}(\mathcal{B})$ constructed in the proof of \cite[Theorem 4.10]{MR3039775}, this local module corresponds to the induction of $X$ to $\mathcal{Z}(\mathcal{B})$.
\item Every simple $Y\in\mathcal{Z}(\mathcal{C})_A^0$ is, as an object of $\mathcal{Z}(\mathcal{C})$, a summand of $I(X)$ for some $X\in\mathcal{O}(\mathcal{B})$.
\end{enumerate}
\end{proposition}

\begin{proof}
\par There is an equivalence of fusion categories $\Z(\mathcal{C})_A\simeq \Z_\mathcal{B}(\mathcal{C})$ by \cite[Lemma 3.5]{MR3039775} because $F_\mathcal{B}$ has the structure of a central surjective functor (see \cite[Definition 2.4]{MR3039775}). Under this equivalence, $\Z(\mathcal{B})\simeq\Z(\mathcal{C})^0_A$ as braided fusion categories as described in the second half of the proof of \cite[Theorem 4.10]{MR3039775}. Moreover, under this equivalence $F_\mathcal{B}$ corresponds to the free module functor $\Z(\mathcal{C})\to\Z(\mathcal{C})_A$, $X\mapsto X\otimes A$ (this is how the equivalence of \cite[Lemma 3.5]{MR3039775} arises), whose right adjoint is the forgetful functor $U\colon\Z(\mathcal{C})_A\to\Z(\mathcal{C})$ sending a module to its underlying object \cite[Theorem 1.6(2)]{MR1936496}. Hence $I_\mathcal{B}$ corresponds to $U$. Since $I(X)=I_\mathcal{B}\widetilde{I}_\mathcal{B}(X)$, both claims follow once we show that for $X\in\mathcal{O}(\mathcal{B})$, the object $\widetilde{I}_\mathcal{B}(X)$ lies in the image of $\mathcal{Z}(\mathcal{B})$ in $\mathcal{Z}_\mathcal{B}(\mathcal{C})$, and in fact corresponds to $I^\mathcal{B}(X)$.

\par To be precise, let $j:\mathcal{Z}(\mathcal{B})\to\mathcal{Z}_\mathcal{B}(\mathcal{C})$ denote the identification of this subcategory, which has a right adjoint $P_j:\mathcal{Z}_\mathcal{B}(\mathcal{C})\to\mathcal{Z}(\mathcal{B})$; as $j$ is a full embedding, $jP_j(W)$ is the largest direct summand of $W$ lying in the image of $j$, for every $W$.  Similarly, we have adjoint functors $\iota:\mathcal{B}\to\mathcal{C}$ and $P_\iota:\mathcal{C}\to\mathcal{B}$.  For simple $X\in\mathcal{C}$, $P_\iota(X)=X$ if $X\in\mathcal{B}$ and $P_\iota(X)=0$ if $X\not\in\mathcal{B}$. Let $F^\mathcal{B}:\mathcal{Z}(\mathcal{B})\to\mathcal{B}$ and $I^\mathcal{B}:\mathcal{B}\to\mathcal{Z}(\mathcal{B})$ be the standard forgetful and induction functors.  In the proof of \cite[Theorem 4.10]{MR3039775}, $j$ is the inclusion of the pairs $(X,\gamma_X)$ with $X\in\mathcal{B}$, i.e.\ it regards an object of $\mathcal{B}$ with a half-braiding relative to $\mathcal{B}$ as an object of $\mathcal{Z}_\mathcal{B}(\mathcal{C})$. Hence $\widetilde{F}_\mathcal{B}\circ j=\iota\circ F^\mathcal{B}$, and taking right adjoints of both sides, $P_j\circ\widetilde{I}_\mathcal{B}\cong I^\mathcal{B}\circ P_\iota$. Therefore if $X\in\mathcal{O}(\mathcal{B})$, then $P_j(\widetilde{I}_\mathcal{B}(X))\cong I^\mathcal{B}(X)$, so $jI^\mathcal{B}(X)$ is a direct summand of $\widetilde{I}_\mathcal{B}(X)$.  Recall that $\mathrm{FPdim}(I^\mathcal{B}(X))=\mathrm{FPdim}(X)\mathrm{FPdim}(\mathcal{B})$.  The same proof as \cite[Proposition 5.4]{MR2183279} applied to the relative center construction shows that
\begin{equation}
    \widetilde{F}_\mathcal{B}(\widetilde{I}_\mathcal{B}(X))\cong\bigoplus_{Y\in\mathcal{O}(\mathcal{B})}Y\otimes X\otimes Y^\ast,
\end{equation}
hence $\mathrm{FPdim}(\widetilde{I}_\mathcal{B}(X))=\mathrm{FPdim}(X)\mathrm{FPdim}(\mathcal{B})=\mathrm{FPdim}(I^\mathcal{B}(X))$.  As $j$ preserves Frobenius--Perron dimensions, the summand $jI^\mathcal{B}(X)$ is all of $\widetilde{I}_\mathcal{B}(X)$, i.e.\ $\widetilde{I}_\mathcal{B}(X)\cong jI^\mathcal{B}(X)$.  Since $I_\mathcal{B}$ corresponds to the forgetful functor, this proves (1).

\par For (2), let $Y\in\mathcal{Z}(\mathcal{B})$ be simple and let $X$ be a simple summand of $F^\mathcal{B}(Y)$, so $X\in\mathcal{O}(\mathcal{B})$.  By adjunction $\mathrm{Hom}(Y,I^\mathcal{B}(X))\cong\mathrm{Hom}(F^\mathcal{B}(Y),X)\neq0$, so $Y$ is a summand of $I^\mathcal{B}(X)$, and $j(Y)$ is a summand of $jI^\mathcal{B}(X)\cong\widetilde{I}_\mathcal{B}(X)$.  Applying $I_\mathcal{B}$, which corresponds to the forgetful functor, shows that the underlying object of the local $A$-module corresponding to $Y$ is a summand of $I(X)$.
\end{proof}


\subsection{Twist bounds in the double}

\par Assume $\mathcal{C}$ is a spherical fusion category. In \cite{MR4655273}, the first author showed that if $g\in\mathcal{O}(\mathcal{C}_\mathrm{pt})$ acting on $\mathcal{O}(\mathcal{C})$ has a fixed-point, the orders of the twists of the simple summands of $I(g)$ divide the order of $g$. Here we sacrifice the strictness of this bound to account for fixed-point free actions by invertible objects.

\begin{proposition}\label{twist}
    Let $\mathcal{C}$ be a spherical fusion category.
    \begin{enumerate}
    \item Let $A$ be a connected \'etale algebra in $\mathcal{Z}(\mathcal{C})$ with $\theta_A=\mathrm{id}_A$, and equip $\mathcal{Z}(\mathcal{C})_A^0$ with the ribbon structure inherited from $\mathcal{Z}(\mathcal{C})$. If $M$ is a simple local $A$-module, then every simple summand of $U(M)$ in $\mathcal{Z}(\mathcal{C})$ has twist $\theta_M$.
    \item Let $\mathcal{B}\subset\mathcal{C}$ be a fusion subcategory and $A\subset I(\mathbbm{1})$ the corresponding connected \'etale algebra, so that $\mathcal{Z}(\mathcal{C})_A^0\simeq\mathcal{Z}(\mathcal{B})$, and let $N\in\mathbb{Z}_{\geq1}$ be the order of the $T$-matrix of $\mathcal{Z}(\mathcal{C})_A^0$ in its inherited ribbon structure. If $X\in\mathcal{O}(\mathcal{B})$ and $Y\subset I(X)$ is simple, then $\theta_Y^N=1$.
    \end{enumerate}
\end{proposition}

\begin{proof}
For the proof of (1), as $\theta_A=\mathrm{id}_A$, the category $\mathcal{Z}(\mathcal{C})_A^0$ inherits a ribbon structure from that of $\mathcal{Z}(\mathcal{C})$ \cite[Theorem 1.17]{MR1936496}, in which the twist of a simple local $A$-module $M$ is the scalar by which $\theta_{U(M)}$ acts; in particular all simple summands of $U(M)$ as an object of $\mathcal{Z}(\mathcal{C})$ have the same twist $\theta_M$ (refer to \cite[Lemma 5.2(b)]{MR3997136}, for example). It remains to prove (2).

Corresponding to the fusion subcategory $\mathcal{B}\subset\mathcal{C}$, there exists a connected \'etale algebra $A\subset I(\mathbbm{1})$ in $\mathcal{Z}(\mathcal{C})$ with $\mathcal{Z}(\mathcal{C})_A^0\simeq\mathcal{Z}(\mathcal{B})$ \cite[Theorem 4.10 and its proof]{MR3039775}. As $A\subset I(\mathbbm{1})$, we have $\theta_A=\mathrm{id}_A$ (Section \ref{subsec:conventions}), so (1) applies. The object $I(X)$ has the structure of a local $A$-module by Proposition \ref{correspond}(1); decomposing this module into simple local modules, every simple $Y\subset I(X)$ is a summand of $U(M)$ for some simple local $A$-module $M$, so $\theta_Y=\theta_M$ by (1), and $\theta_M^N=1$.
\end{proof}

\begin{example}\label{exfix}
    Let $X\in\mathcal{O}(\mathcal{C})$ be invertible generating a cyclic fusion subcategory $\mathcal{B}\subset\mathcal{C}$ with $\mathrm{FPdim}(\mathcal{B})=n\in\mathbb{Z}_{\geq1}$. If $X$ has a fixed-point $Y\in\mathcal{O}(\mathcal{C})$, then the rank one abelian category (i.e., equivalent to $\text{Vec}$ as abelian category) of copies  of $Y$ is a $\mathcal{B}$-module category, hence $\mathcal{B}\simeq\mathrm{Vec}_{C_n}$. We have $T^n=1$ for $\mathcal{Z}(\mathrm{Vec}_{C_n})$ with its canonical ribbon structure, and hence for every ribbon structure: when $n$ is odd the universal grading group $C_n\times C_n$ has odd order, and when $n$ is even the possible sign changes (by $\pm1$-valued characters of the universal grading group) do not affect $T^n$. Therefore $\theta_Z^n=1$ for all simple $Z\subset I(X)$, i.e.\ Proposition \ref{twist} implies \cite[Proposition 3.11]{MR4655273}.
\end{example} 

\subsection{Frobenius--Schur indicators and central objects}\label{subsec:nucentral}

Let $\mathcal{C}$ be a pivotal fusion category with pivotal structure $j$. In this subsection we write the half-braiding of an object $(Z,c)$ of $\mathcal{Z}(\mathcal{C})$ as $c_X\colon Z\otimes X\to X\otimes Z$; in the notation of Section \ref{app:correspond} this is the inverse of $\gamma_Z$. If $V\in\mathcal{O}(\mathcal{C})$ is self-dual and $\varphi\colon V\to V^\ast$ is an isomorphism, then
\begin{equation}\label{eq:nu2}
\varphi^\ast\circ j_V=\nu_2(V)\,\varphi,\qquad \nu_2(V)\in\{\pm1\};
\end{equation}
this is equivalent to the definition of the second Frobenius--Schur indicator in \cite{MR2381536}, transported along $\mathrm{Hom}(\mathbbm{1},V\otimes V)\cong\mathrm{Hom}(V,V^\ast)$. In particular, if $F$ is a pivotal tensor functor and $F(V)$ is simple, then applying $F$ to \eqref{eq:nu2} gives $\nu_2(F(V))=\nu_2(V)$.

\begin{lemma}\label{lem:nucentral}
Let $(Z,c)$ be a self-dual object of $\mathcal{Z}(\mathcal{C})$ whose underlying object $Z$ is simple in $\mathcal{C}$, and let $X\in\mathcal{O}(\mathcal{C})$ be self-dual with $Z\otimes X$ simple. Then $Z\otimes X$ is self-dual and $\nu_2(Z\otimes X)=\nu_2(Z)\,\nu_2(X)$.
\end{lemma}

\begin{proof}
Let $(Z^\ast,c^\vee)$ be the dual of $(Z,c)$ in $\mathcal{Z}(\mathcal{C})$, choose an isomorphism $\varphi_Z\colon(Z,c)\to(Z^\ast,c^\vee)$ in $\mathcal{Z}(\mathcal{C})$, and an isomorphism $\varphi_X\colon X\to X^\ast$ in $\mathcal{C}$. Since $\mathrm{ev}_Z$ is a morphism in $\mathcal{Z}(\mathcal{C})$ \cite[Section 7.13]{MR3242743}, and since $c$ is natural with respect to $\mathrm{ev}_X$ and monoidal, we have for all $Y\in\mathcal{C}$
\begin{align*}
c^\vee_Y&=(\mathrm{ev}_Z\otimes\mathrm{id}_{Y\otimes Z^\ast})\circ(\mathrm{id}_{Z^\ast}\otimes c_Y^{-1}\otimes\mathrm{id}_{Z^\ast})\circ(\mathrm{id}_{Z^\ast\otimes Y}\otimes\mathrm{coev}_Z),\\
c_{X^\ast}^{-1}&=(\mathrm{ev}_X\otimes\mathrm{id}_{Z\otimes X^\ast})\circ(\mathrm{id}_{X^\ast}\otimes c_X\otimes\mathrm{id}_{X^\ast})\circ(\mathrm{id}_{X^\ast\otimes Z}\otimes\mathrm{coev}_X).
\end{align*}
Substituting the second identity into the first with $Y=X^\ast$ gives $c^\vee_{X^\ast}=(c_X)^\ast$ as morphisms $Z^\ast\otimes X^\ast\to X^\ast\otimes Z^\ast$; graphically, both are the rotation of $c_X$ by $\pi$. As $\varphi_Z$ is a morphism in $\mathcal{Z}(\mathcal{C})$ and $c^\vee$ is natural,
\[
(c_X)^\ast\circ(\varphi_Z\otimes\varphi_X)=(\varphi_X\otimes\mathrm{id}_{Z^\ast})\circ c^\vee_X\circ(\varphi_Z\otimes\mathrm{id}_X)=(\varphi_X\otimes\varphi_Z)\circ c_X.
\]
Put $W:=Z\otimes X$ and $\varphi_W:=(\varphi_X\otimes\varphi_Z)\circ c_X\colon W\to X^\ast\otimes Z^\ast=W^\ast$, an isomorphism; so $W$ is self-dual. Since $j$ is monoidal and $(\varphi_X\otimes\varphi_Z)^\ast=\varphi_Z^\ast\otimes\varphi_X^\ast$ under the canonical identifications, \eqref{eq:nu2} for $Z$ and $X$ gives
\[
\varphi_W^\ast\circ j_W=(c_X)^\ast\circ\bigl(\varphi_Z^\ast j_Z\otimes\varphi_X^\ast j_X\bigr)=\nu_2(Z)\nu_2(X)\,(c_X)^\ast\circ(\varphi_Z\otimes\varphi_X)=\nu_2(Z)\nu_2(X)\,\varphi_W .
\]
Now \eqref{eq:nu2} for $W$ gives the claim.
\end{proof}


\section{Stretched near-group doubles}\label{sec:neargroup}

Let $G$ be a finite abelian group and $\mathcal{N}_G$ be a near-group fusion category of type $G+k|G|$ for some $k\in\mathbb{Z}_{\geq1}$. Label the isomorphism classes of invertible objects by $g\in G$ and the isomorphism class of the noninvertible simple objects by $x$. The fusion of $\mathcal{N}_G$ is determined by
\begin{equation}
    xx=k|G|x+\sum_{g\in G}g
\end{equation}
with $\mathrm{FPdim}(\mathcal{N}_G)=|G|(2+kd)$ where $d=\mathrm{FPdim}(x)$ is the largest root of $t^2-k|G|t-|G|$.

\par By \cite[Lemma 3.1]{MR4655273}, $\mathcal{N}_G$ is Galois conjugate to a pseudounitary fusion category. Galois conjugation by an automorphism $\sigma$ of $\mathbb{C}$ maps the double, the algebra $A_G$ and the equivalence $\Lambda$ introduced below to the corresponding data for the conjugate category; it leaves all fusion rules unchanged and replaces each twist $\theta$ by $\sigma(\theta)$, in particular the quadratic form $q$ constructed below by $\sigma\circ q$. We may and do therefore assume throughout Sections \ref{sec:neargroup}--\ref{sec:application} that $\mathcal{N}_G$ is pseudounitary, equipped with its canonical spherical structure, for which $\dim(X)=\mathrm{FPdim}(X)$ for all $X$ \cite[Proposition 9.5.1]{MR3242743}; then $\mathcal{Z}(\mathcal{N}_G)$ is pseudounitary as well. Also note that $k^2|G|^2<k^2|G|^2+4|G|<(k|G|+2)^2$, and $k^2|G|^2+4|G|=(k|G|+1)^2$ would force $4|G|=2k|G|+1$, which has no solutions by parity. Hence the discriminant $k^2|G|^2+4|G|$ is not a square, and $d\notin\mathbb{Q}$.


\subsection{Induction via local modules}\label{subsec:localpt}

Throughout this section we put $n:=|G|$ and $\alpha:=2+kd$, and in this subsection $k\geq1$ is arbitrary. We write $\mu_m$ for the group of $m$-th roots of unity. From $d^2=knd+n$ we get the identities
\begin{equation}\label{eq:ids}
  d^{2}=n(\alpha-1),\qquad
  \alpha^{2}=(4+k^2n)(\alpha-1),\qquad\text{and}\qquad
  1+(\alpha-1)^{2}=(2+k^{2}n)(\alpha-1).
\end{equation}

\par Let $\mathcal{B}:=(\mathcal{N}_G)_\mathrm{pt}$. Every $g\in\mathcal{O}(\mathcal{B})$ fixes $x$, so the full subcategory of $\mathcal{N}_G$ consisting of direct sums of copies of $x$ is a $\mathcal{B}$-module category of rank one, i.e.\ there is a fiber functor $\mathcal{B}\to\mathrm{Vec}$. Hence $\mathcal{B}\simeq\mathrm{Vec}_G$ (compare Example \ref{exfix}). Let $A_G\subset I(\mathbbm{1})$ be the connected \'etale algebra in $\mathcal{Z}(\mathcal{N}_{G})$ corresponding to $\mathcal{B}$ \cite[Theorem 4.10]{MR3039775}. It has $\dim(A_G)=\dim(\mathcal{N}_G)/\dim(\mathcal{B})=\alpha$, and there is a braided equivalence
\begin{equation}\label{theeq}
\Lambda\colon\mathcal{Z}(\mathrm{Vec}_G)\xrightarrow{\ \sim\ }\mathcal{Z}(\mathcal{N}_G)_{A_G}^0
\end{equation}
(see Section \ref{app:correspond}). We distinguish two forgetful functors. As in Section \ref{app:correspond}, the functor $U\colon\mathcal{Z}(\mathcal{N}_G)_{A_G}\to\mathcal{Z}(\mathcal{N}_G)$ sends an $A_G$-module $M$ to its underlying object $U(M)\in\mathcal{Z}(\mathcal{N}_G)$. The functor $F\colon\mathcal{Z}(\mathcal{N}_G)\to\mathcal{N}_G$ is the forgetful functor to $\mathcal{N}_G$, and we call $F(Y)\in\mathcal{N}_G$ the \emph{forgetful image} of $Y\in\mathcal{Z}(\mathcal{N}_G)$. Since $A_G\subset I(\mathbbm{1})$, we have $\theta_{A_G}=\mathrm{id}$ (Section \ref{subsec:conventions}). Hence $\mathcal{Z}(\mathcal{N}_G)^0_{A_G}$ is a ribbon category with twists inherited from $\mathcal{Z}(\mathcal{N}_G)$ \cite[Theorem 1.17]{MR1936496}, and the dimension of a local module $M$ is $\dim(U(M))/\alpha=\dim(F(U(M)))/\alpha$ \cite[Theorem 1.18]{MR1936496}. The simple objects of $\mathcal{Z}(\mathrm{Vec}_G)$ are the pairs $\gamma=(g,\chi)\in\Gamma:=G\times\widehat{G}$, consisting of an invertible object $g$ and a character $\chi$ recording its half-braiding. They satisfy $(g,\chi)\otimes(h,\psi)=(gh,\chi\psi)$ and $(g,\chi)^\ast=(g^{-1},\chi^{-1})$. With its canonical ribbon structure, $\mathcal{Z}(\mathrm{Vec}_G)$ has all dimensions equal to $1$ and twists $\theta_{(g,\chi)}=\chi(g)$. In the convention \eqref{eq:balancing}, its $S$-matrix is
\begin{equation}\label{eq:S0}
S^0_{(g,\chi),(h,\psi)}=\overline{\chi(h)\psi(g)}.
\end{equation}
The equivalence $\Lambda$ in \eqref{theeq} is one of ribbon fusion categories, where $\mathcal{Z}(\mathrm{Vec}_G)$ carries its canonical ribbon structure and $\mathcal{Z}(\mathcal{N}_G)^0_{A_G}$ the one inherited from $\mathcal{Z}(\mathcal{N}_G)$. Indeed, transported along $\Lambda$, the two differ by a $\{\pm1\}$-valued character of the universal grading group, and this character multiplies the dimension of each simple object by its value on the degree of that object. Both the canonical dimensions and the inherited dimensions $\dim(U(M))/\alpha$ are positive, so the character is trivial.

\par For $\gamma\in\Gamma$ we write $M_\gamma:=\Lambda(\gamma)$ for the corresponding simple local $A_G$-module and $U_\gamma:=U(M_\gamma)$, an object of $\mathcal{Z}(\mathcal{N}_G)$ of dimension $\alpha$, whose forgetful image $F(U_\gamma)$ is an object of $\mathcal{N}_G$. The functors involved fit into the diagram
\begin{equation}\label{eq:functors}
\begin{tikzcd}[column sep=large, row sep=large]
\mathrm{Vec}_G \arrow[r, "I^{\mathrm{Vec}_G}"] \arrow[d, hook, "\iota"'] & \mathcal{Z}(\mathrm{Vec}_G) \arrow[r, "\Lambda", "\simeq"'] & \mathcal{Z}(\mathcal{N}_G)^0_{A_G} \arrow[d, "U"] \\
\mathcal{N}_G \arrow[rr, "I"'] & & \mathcal{Z}(\mathcal{N}_G)
\end{tikzcd}
\end{equation}
where $\iota$ is the inclusion of $\mathcal{B}\simeq\mathrm{Vec}_G$ and $I^{\mathrm{Vec}_G}$ is the induction to $\mathcal{Z}(\mathrm{Vec}_G)$, so that $I^{\mathrm{Vec}_G}(g)=\bigoplus_{\chi\in\widehat{G}}(g,\chi)$. By Proposition \ref{correspond}(1), for every $g\in G$ we have $I(g)\cong U\Lambda I^{\mathrm{Vec}_G}(g)=\bigoplus_{\chi}U_{(g,\chi)}$; that is, the diagram \eqref{eq:functors} commutes up to isomorphism on objects of $\mathrm{Vec}_G$.  The structure of $\mathcal{Z}(\mathcal{N}_G)$ is dictated by the following.

\begin{proposition}\label{prop:local}
Let $g\in G$.
\begin{enumerate}
\item $I(g)\cong\bigoplus_{\chi\in\widehat{G}}U_{(g,\chi)}$, and $I(g)$ is multiplicity free.
\item There is a unique $t_g\in\widehat{G}$ such that $U_{(g,t_g)}$ is not simple.
\begin{itemize}
\item It decomposes as $U_{(g,t_g)}\cong Y_g\oplus Y_g'$ with $Y_g\not\cong Y_g'$ simple.
\item Here $F(Y_g)=g+\ell_{1,g}x$ and $F(Y_g')=g+\ell_{2,g}x$ for some $\ell_{1,g},\ell_{2,g}\in\mathbb{Z}_{\geq0}$ with $\ell_{1,g}+\ell_{2,g}=k$.
\item For $g=e$ we have $t_e=1$ and $U_{(e,1)}=A_G$. We take $Y_e:=\mathbbm{1}$, so that $\ell_{1,e}=0$ and $A_G\cong\mathbbm{1}\oplus Y_e'$.
\end{itemize}
\item Put $\mathcal{T}:=\{(h,t_h):h\in G\}\subset\Gamma$. For $\gamma\in\Gamma\setminus\mathcal{T}$ the object $Y_\gamma:=U_\gamma$ is simple, and there is a fixed-point-free involution $\gamma\mapsto\gamma^\sharp$ of $\Gamma\setminus\mathcal{T}$ such that, for $\gamma,\gamma'\in\Gamma\setminus\mathcal{T}$,
\[
U_\gamma\cong U_{\gamma'}\iff\gamma'\in\{\gamma,\gamma^\sharp\}\iff F(U_\gamma)\cong F(U_{\gamma'}).
\]
If $\chi\neq t_g$ and $(g,\chi)^\sharp=(h,\psi)$, then $h\neq g$ and $F(Y_{(g,\chi)})=g+h+kx$; moreover, the assignment $\chi\mapsto h$ is a bijection $\widehat{G}\setminus\{t_g\}\to G\setminus\{g\}$.
\item Every simple summand of $U_{(g,\chi)}$ has twist $\chi(g)$. In particular $\theta_{Y_g}=\theta_{Y_g'}=t_g(g)$, and $\theta_{Y_{(g,\chi)}}=\chi(g)$ for $\chi\neq t_g$.
\item Let $W$ be a simple object of $\mathcal{Z}(\mathcal{N}_G)$ which is not a summand of $\bigoplus_{f\in G}I(f)$. Then $F(W)=m_Wx$ for some $m_W\in\mathbb{Z}_{\geq1}$.
\item $U_\gamma^\ast\cong U_{\gamma^\ast}$ for every $\gamma\in\Gamma$. In particular $Y_{(g,\chi)}^\ast\cong Y_{(g^{-1},\chi^{-1})}$ for $\chi\neq t_g$.
\end{enumerate}
\end{proposition}

\begin{proof}
The isomorphism $I(g)\cong\bigoplus_\chi U_{(g,\chi)}$ was noted after \eqref{eq:functors}. By \cite[Proposition 5.4]{MR2183279},
\begin{equation}\label{eq:FIg}
F(I(g)) = F\Bigl(\bigoplus_\chi U_{(g,\chi)}\Bigr)=\sum_{h\in G}hgh^{-1}+xgx^\ast=ng+\sum_{h\in G}h+knx,
\end{equation}
so $[F(I(g)),g]=n+1$ and $[F(I(g)),h]=1$ for $h\neq g$. By Frobenius reciprocity, a simple object $Y$ of $\mathcal{Z}(\mathcal{N}_G)$ occurs in $I(g)$ with multiplicity $[Y,I(g)]=[F(Y),g]$, hence
\begin{equation}\label{eq:musum}
\sum_{Y \in \mathcal{O}(\mathcal{Z}(\mathcal{N}_G))}[F(Y),g]^2=[F(I(g)),g]=n+1.
\end{equation}

\par Since $\dim(U_{(g,\chi)})=2+kd$ and $d\notin\mathbb{Q}$, the object $F(U_{(g,\chi)})$ is the sum of $kx$ and exactly two invertible objects, counted with multiplicity. Every simple summand $Y$ of $U_{(g,\chi)}$ is a summand of $I(g)$, so $[F(Y),g]\geq1$. Comparing with \eqref{eq:FIg}, exactly one character, which we call $t_g$, has $[F(U_{(g,t_g)}),g]=2$, since the $n$ integers $[F(U_{(g,\chi)}),g]\in\{1,2\}$ sum to $n+1$; thus $F(U_{(g,t_g)})=2g+kx$. For $\chi\neq t_g$ we have $F(U_{(g,\chi)})=g+h+kx$ with $h\neq g$, and $\chi\mapsto h$ is a bijection $\widehat{G}\setminus\{t_g\}\to G\setminus\{g\}$. Moreover, for $\chi\neq t_g$ the object $U_{(g,\chi)}$ is simple and occurs in $I(g)$ with multiplicity one, since each of its simple summands has $g$ in its forgetful image and $g$ occurs there only once. The summands $U_{(g,\chi)}$ with $\chi\neq t_g$ therefore contribute $n-1$ to \eqref{eq:musum}, so the simple summands of $U_{(g,t_g)}$ contribute $2$. We rule out the two alternatives. If $U_{(g,t_g)}$ were simple, it would contribute $[F(U_{(g,t_g)}),g]^2=4$; and if $U_{(g,t_g)}\cong Y\oplus Y$, then $Y$ would occur in $I(g)$ with multiplicity at least $2$, whereas $[F(Y),g]=1$. Hence $U_{(g,t_g)}\cong Y_g\oplus Y_g'$ with $Y_g\not\cong Y_g'$, each occurring in $I(g)$ with multiplicity one, and $I(g)$ is multiplicity free. As $F(U_{(g,t_g)})=2g+kx$ and each of $F(Y_g)$, $F(Y'_g)$ contains $g$ once, we may write $F(Y_g)=g+\ell_{1,g}x$ and $F(Y'_g)=g+\ell_{2,g}x$ with $\ell_{1,g}+\ell_{2,g}=k$. This proves (1) and the first two parts of (2).

\par The unit object of $\mathcal{Z}(\mathcal{N}_G)_{A_G}^0$ is $M_{(e,1)}$, whose image under $U$ is $A_G=U_{(e,1)}$. Since $A_G$ contains $\mathbbm{1}$ and $\alpha\neq1$, the object $U_{(e,1)}=A_G$ is not simple. As $U_{(e,\chi)}$ is simple for every $\chi\neq t_e$, this forces $t_e=1$. We may therefore take $Y_e:=\mathbbm{1}$ and write $A_G=\mathbbm{1}\oplus Y'_e$, which gives the claims about $g=e$.

\par For $\gamma\in\Gamma\setminus\mathcal{T}$ the object $Y_\gamma=U_\gamma$ is simple by the above. Let $\chi\neq t_g$ and put $Y:=Y_{(g,\chi)}$, with $F(Y)=g+h+kx$. Since $[F(Y),h]=1$, $Y$ occurs in $I(h)$ with multiplicity one, hence is a summand of exactly one $U_{(h,\psi)}$ by (1). We have $\psi\neq t_h$, since the forgetful images of $Y_h$ and $Y_h'$ do not contain $g$; hence $Y\cong Y_{(h,\psi)}$. As $[F(Y),f]=0$ for $f\notin\{g,h\}$, $Y$ is not a summand of $I(f)$ for any other $f$. Consequently $(g,\chi)$ and $(h,\psi)$ are the only labels in $\Gamma\setminus\mathcal{T}$ of the object $Y$, and putting $(g,\chi)^\sharp:=(h,\psi)$ defines an involution of $\Gamma\setminus\mathcal{T}$ with $U_\gamma\cong U_{\gamma'}$ if and only if $\gamma'\in\{\gamma,\gamma^\sharp\}$. It changes the first coordinate, so it has no fixed points. Finally, $Y$ is the unique simple summand of $I(g)$ whose forgetful image contains $h$, by the bijection $\chi\mapsto h$; so if $F(U_{\gamma'})\cong F(Y)$ then $U_{\gamma'}$, which is a summand of $I(g)$ by Frobenius reciprocity, is isomorphic to $Y$. This proves (3).

\par By Proposition \ref{twist}(1), every simple summand of $U_{(g,\chi)}$ has the twist of $M_{(g,\chi)}$ in the inherited ribbon structure, which is the canonical one; so it has twist $\theta_{(g,\chi)}=\chi(g)$, which proves (4). If a simple object $W$ has $[F(W),f]\neq0$ for some $f\in G$, then $W$ is a summand of $I(f)$ by Frobenius reciprocity; so for $W$ as in (5), $F(W)$ contains no invertible object, that is, $F(W)=m_Wx$ with $m_W\geq1$. Finally, the dual of $M_\gamma$ is $M_{\gamma^\ast}$, and its underlying object is $U_\gamma^\ast$ (see the proof of \cite[Theorem 1.15]{MR1936496}); as duality preserves simplicity, $\chi^{-1}\neq t_{g^{-1}}$ whenever $\chi\neq t_g$, and this gives (6).
\end{proof}

The $S$-matrix of $\mathcal{Z}(\mathcal{N}_G)$ is related to that of $\mathcal{Z}(\mathrm{Vec}_G)$ by the following branching rule.

\begin{theorem}[{\cite[Theorem 4.1]{MR1936496}}]\label{thm:branch}
Let $\mathcal{A}$ be a modular fusion category and $A\in\mathcal{A}$ a connected \'etale algebra with $\theta_A=\mathrm{id}_A$ and $\dim(A)\neq0$. For $X\in\mathcal{O}(\mathcal{A})$ and a simple local $A$-module $M$,
\begin{equation}\label{eq:branch}
S_{X,U(M)}=\dim(A)\sum_{M'}[X,U(M')]\,S^{0}_{M',M},
\end{equation}
where $M'$ runs over the simple local $A$-modules and $S^0$ is the $S$-matrix of $\mathcal{A}^0_A$ with its inherited ribbon structure.
\end{theorem}

\begin{proof}
Let $\Phi\colon K_0(\mathcal{A}_A^0)\to K_0(\mathcal{A})$ be induced by the forgetful functor, so that $\Phi[M]=\sum_X[X,U(M)]\,[X]$, and let $\tilde s$ and $\tilde s^A$ be the operators $[V]\mapsto\sum_WS_{V,W}[W]$ on the two Grothendieck rings \cite[(4.2)]{MR1936496}. The hypotheses of \cite{MR1936496} (a rigid $\mathcal{A}$-algebra with trivial twist, in the terminology there) hold here, since the pairing $A\otimes A\to A\to\mathbbm{1}$ is nondegenerate for connected \'etale $A$ \cite[Remark 3.4]{MR3039775}. So \cite[Theorem 4.1]{MR1936496} applies and states that $\dim(A)\,\Phi\circ\tilde s^A=\tilde s\circ\Phi$. Applying both sides to $[M]$, comparing the coefficients of $[X]$, and using the symmetry of both $S$-matrices gives \eqref{eq:branch}. It remains to compare conventions. The Hopf-link convention of \cite{MR1936496} agrees with \eqref{eq:balancing} up to replacing $S_{V,W}$ by $S_{V^\ast,W}$ in both $\mathcal{A}$ and $\mathcal{A}^0_A$. Since $U(M'^\ast)\cong U(M')^\ast$ (see the proof of \cite[Theorem 1.15]{MR1936496}), we have $[X^\ast,U(M')]=[X,U(M'^\ast)]$, and $M'\mapsto M'^\ast$ permutes the simple local modules. Hence \eqref{eq:branch} is invariant under this replacement, so it also holds in the convention \eqref{eq:balancing}.
\end{proof}

\begin{corollary}\label{cor:branch}
Let $\gamma\in\Gamma$ and $g,h\in G$.
\begin{enumerate}
\item Let $W$ be a simple object which is not a summand of $\bigoplus_{f\in G}I(f)$. Then $S_{W,U_\gamma}=0$. In particular $S_{W,Y_g}=-S_{W,Y_g'}$, and $S_{W,Y_{(g,\chi)}}=0$ for $\chi\neq t_g$.
\item For $X\in\{Y_h,Y_h'\}$ we have $S_{X,U_\gamma}=\alpha\,S^0_{(h,t_h),\gamma}$. In particular, for $\chi\neq t_g$,
\[
S_{Y_h,Y_{(g,\chi)}}=S_{Y_h',Y_{(g,\chi)}}=\alpha\,S^0_{(h,t_h),(g,\chi)}=\alpha\,S^0_{(h,t_h),(g,\chi)^\sharp}.
\]
\item We have $S_{Y_g,Y_h}=S_{Y_g',Y_h'}$ and $S_{Y_g,Y_h'}=S_{Y_g',Y_h}$, and
\[
S_{Y_g,Y_h}+S_{Y_g,Y_h'}=\alpha\,S^0_{(g,t_g),(h,t_h)}.
\]
\end{enumerate}
\end{corollary}

\begin{proof}
We apply Theorem \ref{thm:branch} with $\mathcal{A}=\mathcal{Z}(\mathcal{N}_G)$ and $A=A_G$, which satisfies its hypotheses since $\theta_{A_G}=\mathrm{id}$ and $\dim(A_G)=\alpha>0$. If $W$ is a simple object which is not a summand of $\bigoplus_{f\in G}I(f)$, then $[W,U(M')]=0$ for every simple local $A_G$-module $M'$ by Proposition \ref{prop:local}(1) and (5), so \eqref{eq:branch} gives $S_{W,U_\gamma}=0$. As $U_{(g,t_g)}\cong Y_g\oplus Y_g'$ and $U_{(g,\chi)}=Y_{(g,\chi)}$ for $\chi\neq t_g$ by Proposition \ref{prop:local}(2) and (3), this gives the remaining identities in (1). If instead $X\in\{Y_h,Y_h'\}$, then $[X,U(M')]=\delta_{M',M_{(h,t_h)}}$ by Proposition \ref{prop:local}(1)--(3), and \eqref{eq:branch} gives $S_{X,U_\gamma}=\alpha\,S^0_{(h,t_h),\gamma}$; the identities displayed in (2) follow, as $U_{(g,\chi)}=Y_{(g,\chi)}=U_{(g,\chi)^\sharp}$ for $\chi\neq t_g$. Finally, put $s_1=S_{Y_g,Y_h}$, $s_2=S_{Y_g,Y_h'}$, $s_3=S_{Y_g',Y_h}$ and $s_4=S_{Y_g',Y_h'}$. Applying (2) to $X\in\{Y_g,Y_g'\}$ with $\gamma=(h,t_h)$ and to $X\in\{Y_h,Y_h'\}$ with $\gamma=(g,t_g)$, and using the symmetry of $S$ and $S^0$, we get $s_1+s_2=s_3+s_4=s_1+s_3=s_2+s_4=\alpha\,S^0_{(g,t_g),(h,t_h)}$. Hence $s_1=s_4$ and $s_2=s_3$, which proves (3).
\end{proof}

We now construct the quadratic form $q$ on $G$ for every $k$.

\begin{proposition}\label{prop:tq}
\leavevmode
\begin{enumerate}
\item If $(g,\chi)^\sharp=(h,\psi)$, then $\psi(f)=\chi(f)\,t_f(gh^{-1})$ for all $f\in G$, and $\chi(g)=\psi(h)$.
\item The map $t\colon G\to\widehat{G}$, $g\mapsto t_g$, is a group homomorphism.
\item The function $q(g):=\theta_{Y_g}=\theta_{Y_g'}=t_g(g)$ is a quadratic form on $G$. Its associated bicharacter is
\begin{equation}\label{eq:bt}
b(g,h):=q(gh)q(g)^{-1}q(h)^{-1}=t_g(h)\,t_h(g),
\end{equation}
and $S^0_{(g,t_g),(h,t_h)}=\overline{b(g,h)}$.
\end{enumerate}
\end{proposition}

\begin{proof}
Suppose that $(g,\chi)^\sharp=(h,\psi)$. Corollary \ref{cor:branch}(2) with $X=Y_f$ gives $S^0_{(f,t_f),(g,\chi)}=S^0_{(f,t_f),(h,\psi)}$, which by \eqref{eq:S0} says that $t_f(g)\chi(f)=t_f(h)\psi(f)$, that is, $\psi(f)=\chi(f)\,t_f(gh^{-1})$; moreover the simple object $Y_{(g,\chi)}\cong Y_{(h,\psi)}$ has twist $\chi(g)=\psi(h)$ by Proposition \ref{prop:local}(4). This proves (1). To see that $t$ is a homomorphism, note first that $t_e=1$ by Proposition \ref{prop:local}(2). If $u\in G\setminus\{e\}$, then by Proposition \ref{prop:local}(3) there is $\chi\neq t_u$ with $(u,\chi)^\sharp=(e,\psi)$ for some $\psi$, and by (1) the function $f\mapsto t_f(u)=\psi(f)\chi(f)^{-1}$ is a character of $G$. Hence $t_{ff'}(u)=t_f(u)\,t_{f'}(u)$ for all $f,f',u\in G$, which proves (2). Finally, $q(g)=t_g(g)$ by Proposition \ref{prop:local}(4), so (2) gives $q(g^{-1})=t_{g^{-1}}(g^{-1})=t_g(g)=q(g)$ and $q(gh)=t_{gh}(gh)=t_g(g)\,t_g(h)\,t_h(g)\,t_h(h)$. The latter identity is \eqref{eq:bt}, and its right-hand side $t_g(h)\,t_h(g)$ is a symmetric bicharacter, so $q$ is a quadratic form. The identity $S^0_{(g,t_g),(h,t_h)}=\overline{b(g,h)}$ is then \eqref{eq:S0}.
\end{proof}

\begin{proposition}\label{prop:block}
Let $g\in G$ have order $m$. Then:
\begin{enumerate}
\item $\theta_Y^{\,m}=1$ for every simple $Y\subseteq I(g)$; in particular, all values of $q$ and $b$ lie in $\mu_{\exp(G)}$;
\item for every $\zeta\in\mu_m$,
\begin{equation}\label{eq:fibre}
  \#\bigl\{\chi\in\widehat{G}\setminus\{t_g\}\ :\ \theta_{Y_{(g,\chi)}}=\zeta\bigr\}
  \ =\ \frac nm-\delta_{\zeta,\,q(g)} .
\end{equation}
\end{enumerate}
\end{proposition}

\begin{proof}
By Proposition \ref{prop:local}(4), the twists of the simple summands of $I(g)$ are the values $\chi(g)$ for $\chi\in\widehat{G}$: the character $t_g$ contributes $Y_g$ and $Y_g'$, both of twist $t_g(g)=q(g)$, and each $\chi\neq t_g$ contributes the single object $Y_{(g,\chi)}$, of twist $\chi(g)$. The map $\chi\mapsto\chi(g)$ is a surjective homomorphism $\widehat{G}\to\mu_m$, each of whose fibres has $n/m$ elements. Both claims follow.
\end{proof}

\par Since $F$ is a tensor functor, it restricts to a group homomorphism $\mathcal{O}(\mathcal{Z}(\mathcal{N}_G)_\mathrm{pt})\to G$. We denote its image by $K$.

\begin{lemma}\label{lem:liftable}
\leavevmode
\begin{enumerate}
\item The homomorphism $\mathcal{O}(\mathcal{Z}(\mathcal{N}_G)_\mathrm{pt})\to G$ is injective. Moreover, $g\in K$ if and only if $\min(\ell_{1,g},\ell_{2,g})=0$. For $g\in K$ we label so that $\ell_{1,g}=0$; then $Y_g$ is the invertible object over $g$.
\item For $g\in K$ we have $Y_g\otimes Y_e'\cong Y_g'$. Moreover, the free $A_G$-module $A_G\otimes Y_g$ is isomorphic to $M_{(g,t_g)}$ in $\mathcal{Z}(\mathcal{N}_G)_{A_G}$; in particular $A_G\otimes Y_g\cong U_{(g,t_g)}$ in $\mathcal{Z}(\mathcal{N}_G)$.
\item For $u\in K$ and $\gamma\in\Gamma$ we have $Y_u\otimes U_\gamma\cong U_{(u,t_u)\gamma}$. In particular $Y_u\otimes Y_{(g,\chi)}\cong Y_{(ug,t_u\chi)}$ for $\chi\neq t_g$.
\item Let $u\in K$ and let $Y$ be a simple summand of $U_{(h,\psi)}$. Then the double braiding $c_{Y,Y_u}c_{Y_u,Y}$ is the scalar $t_u(h)\,\psi(u)$.
\end{enumerate}
\end{lemma}

\begin{proof}
Let $Y$ be invertible with $F(Y)=g$. Then $Y$ is a summand of $I(g)$ of dimension $1$, so by Proposition \ref{prop:local} it is $Y_g$ with $\ell_{1,g}=0$ or $Y_g'$ with $\ell_{2,g}=0$; conversely, a summand of $I(g)$ with forgetful image $g$ has dimension $1$ and is invertible. Hence $g\in K$ if and only if $\min(\ell_{1,g},\ell_{2,g})=0$, and since for $g=e$ the only such object is $\mathbbm{1}$, the homomorphism $\mathcal{O}(\mathcal{Z}(\mathcal{N}_G)_\mathrm{pt})\to G$ is injective. This proves (1). Now let $g\in K$. The object $Y_g\otimes Y_e'$ is simple with $F(Y_g\otimes Y_e')=g+kx$, so it is a summand of $I(g)$ of dimension $1+kd$, and hence it is $Y_g'$ by Proposition \ref{prop:local}. Consider the free $A_G$-module $A_G\otimes Y_g$, where the tensor product is taken in $\mathcal{Z}(\mathcal{N}_G)$ and $A_G$ acts through its multiplication \cite[Theorem 1.6]{MR1936496}. Since $A_G\cong\mathbbm{1}\oplus Y'_e$, its underlying object is $U(A_G\otimes Y_g)\cong Y_g\oplus Y_e'\otimes Y_g\cong U_{(g,t_g)}$, and $\mathrm{Hom}_{A_G}(A_G\otimes Y_g,M_{(g,t_g)})\cong\mathrm{Hom}(Y_g,U_{(g,t_g)})\neq0$ by \cite[Theorem 1.6(2)]{MR1936496}. As $\mathcal{Z}(\mathcal{N}_G)_{A_G}$ is semisimple and $M_{(g,t_g)}$ is simple, a nonzero morphism $A_G\otimes Y_g\to M_{(g,t_g)}$ exhibits $M_{(g,t_g)}$ as a direct summand of $A_G\otimes Y_g$, and its complement has zero underlying object because $U(A_G\otimes Y_g)\cong U_{(g,t_g)}=U(M_{(g,t_g)})$. Hence $A_G\otimes Y_g\cong M_{(g,t_g)}$, which proves (2). Next let $u\in K$ and $\gamma\in\Gamma$, and write $\otimes_{A_G}$ for the tensor product of $\mathcal{Z}(\mathcal{N}_G)_{A_G}$. By \cite[Theorem 1.6(4)]{MR1936496} we have $U\bigl((A_G\otimes Y_u)\otimes_{A_G}M_\gamma\bigr)\cong Y_u\otimes U_\gamma$, while (2) and the fact that $\Lambda$ is a tensor functor give $(A_G\otimes Y_u)\otimes_{A_G}M_\gamma\cong M_{(u,t_u)}\otimes_{A_G}M_\gamma\cong M_{(u,t_u)\gamma}$. If $\gamma=(g,\chi)$ with $\chi\neq t_g$, then $t_u\chi\neq t_ut_g=t_{ug}$ by Proposition \ref{prop:tq}(2), so $U_{(ug,t_u\chi)}=Y_{(ug,t_u\chi)}$; this proves (3). Finally, let $Y$ be a simple summand of $U_{(h,\psi)}$. By (3), $Y_u\otimes Y$ is a simple summand of $U_{(uh,t_u\psi)}$, so by Proposition \ref{prop:local}(4) its twist is $(t_u\psi)(uh)=t_u(u)\,t_u(h)\,\psi(u)\,\psi(h)$. As $\theta_{Y_u}=t_u(u)$ and $\theta_Y=\psi(h)$, the balancing axiom $\theta_{Y_u\otimes Y}=(\theta_{Y_u}\otimes\theta_Y)\,c_{Y,Y_u}c_{Y_u,Y}$ (see e.g.\ \cite[Section 8.10]{MR3242743}) gives (4).
\end{proof}

Several results below need every $g\in G$ to lift to an invertible object of the double. We isolate this as a hypothesis.

\begin{definition}\label{hyp:H}
We say that $\mathcal{N}_G$ satisfies \emph{hypothesis} \eqref{eq:star} if $K=G$, that is,
\begin{equation}\label{eq:star}\tag{H}
  F\ \text{restricts to a group isomorphism}\
  \mathcal{O}(\mathcal{Z}(\mathcal{N}_G)_\mathrm{pt})\xrightarrow{\ \sim\ }G .
\end{equation}
In that case we label so that $\ell_{1,g}=0$ for all $g\in G$ (Lemma \ref{lem:liftable}(1)). Then the objects $Y_g$ are invertible and form a group isomorphic to $G$ via $F$, and Lemma \ref{lem:liftable}(2)--(4) applies to every element of $G$.
\end{definition}


\subsection{The quadratic form for cyclic \texorpdfstring{$G$}{G}}\label{subsec:prim}

Throughout this subsection $k\geq1$ is arbitrary. The arguments use only Subsection \ref{subsec:localpt}; in particular they do not assume that $\mathcal{Z}(\mathcal{N}_G)$ has nontrivial invertible objects. For any finite abelian group $G$ we write $\operatorname{rad}(b):=\{g\in G:b(g,-)\equiv1\}$. When $G$ is cyclic, $t$, $q$ and $b$ are determined completely.

\begin{proposition}\label{cor:cyclic}
Let $G=C_n$, identified with $\mathbb{Z}/n\mathbb{Z}$ via a generator. Then $\kappa:=q(1)$ is a primitive $n$-th root of unity, the homomorphism $t\colon G\to\widehat{G}$ of Proposition \ref{prop:tq}(2) is an isomorphism, and
\begin{equation}\label{eq:qbcyclic}
  t_i(j)=\kappa^{\,ij},\qquad q(j)=\kappa^{\,j^{2}},\qquad b(i,j)=\kappa^{\,2ij},\qquad
  \operatorname{rad}(b)=\begin{cases}\{0\}, & n \text{ odd},\\ \{0,\tfrac n2\}, & n\text{ even}.\end{cases}
\end{equation}
For $n$ even, $q(\tfrac n2)=1$ if $4\mid n$ and $q(\tfrac n2)=-1$ if
$n\equiv2\pmod4$.
\end{proposition}
 
\begin{proof}
Let $g\in G$ correspond to $1$, so that $\kappa=q(g)$, and let $p$ be a prime dividing $n$. By Proposition \ref{prop:block}(1), $\kappa\in\mu_n$, and we show that $\kappa\notin\mu_{n/p}$. Let $L\subset G$ be the subgroup of order $n/p$. For $\chi\in\widehat{G}\setminus\{t_g\}$ write $(g,\chi)^\sharp=(h_\chi,\psi_\chi)$. By Proposition \ref{prop:local}(3), $\chi\mapsto h_\chi$ is a bijection onto $G\setminus\{g\}$, which contains $L$. If $h_\chi\in L$, then $Y_{(g,\chi)}\cong Y_{(h_\chi,\psi_\chi)}$ is a summand of $I(h_\chi)$, so $\theta_{Y_{(g,\chi)}}\in\mu_{n/p}$ by Proposition \ref{prop:block}(1). Therefore, if $\kappa\in\mu_{n/p}$, summing \eqref{eq:fibre} over $\zeta\in\mu_{n/p}$ gives
\[
  \frac np=|L|\ \le\ \#\bigl\{\chi\neq t_g:\theta_{Y_{(g,\chi)}}\in\mu_{n/p}\bigr\}
  \ =\ \sum_{\zeta\in\mu_{n/p}}\bigl(1-\delta_{\zeta,\kappa}\bigr)\ =\ \frac np-1,
\]
a contradiction. Hence $\kappa=t_1(1)$ has order $n$, so the character $t_1$ has order $n$ and generates $\widehat{G}\cong C_n$. As $t$ is a homomorphism by Proposition \ref{prop:tq}(2), it is surjective, hence an isomorphism. Moreover $t_i=t_1^{\,i}$ and $t_1(j)=t_1(1)^j=\kappa^j$, so $t_i(j)=\kappa^{ij}$; then $q(j)=t_j(j)=\kappa^{j^2}$, and $b(i,j)=t_i(j)t_j(i)=\kappa^{2ij}$ by \eqref{eq:bt}. Thus $i\in\operatorname{rad}(b)$ if and only if $n\mid 2ij$ for all $j$, i.e.\ $n\mid 2i$.
 
Finally $q(\tfrac n2)=\kappa^{n^2/4}$. If $4\mid n$, then $n\mid n^2/4$, so $q(\tfrac n2)=1$. If $n=2n'$ with $n'$ odd, then $n\mid n'^2-n'$, so $q(\tfrac n2)=\kappa^{n'}=-1$.
\end{proof}

\subsection{Frobenius--Schur indicators and twist identities}\label{subsec:indicators}

We compute $\nu_r(x)$ for $r$ prime to the Frobenius--Schur exponent, and derive identities for the twists of the simple objects of $\mathcal{Z}(\mathcal{N}_G)$ that do not lie over invertible objects. Throughout this subsection $k\geq1$ is arbitrary, and \eqref{eq:star} is not assumed. Put
\[
\Delta:=k^2n^2+4n=n(k^2n+4),
\]
so that $d=\frac12(kn+\sqrt\Delta)$, and let $d':=\frac12(kn-\sqrt\Delta)$ and $\alpha':=2+kd'$ be the Galois conjugates of $d$ and $\alpha$. Then
\begin{equation}\label{eq:ddprime}
  d+d'=kn,\qquad dd'=-n,\qquad d-d'=\sqrt\Delta,\qquad \alpha\alpha'=4+k^2n .
\end{equation}
Let $N:=\mathrm{FSexp}(\mathcal{N}_G)$ be the order of the $T$-matrix of $\mathcal{Z}(\mathcal{N}_G)$, and let $\mathfrak{d}$ be the discriminant of $\mathbb{Q}(d)=\mathbb{Q}(\sqrt\Delta)$; its absolute value is the conductor $\mathfrak{f}$ of $\mathbb{Q}(\sqrt\Delta)$. For $r\in\mathbb{Z}_{\geq1}$ write $G[r]:=\{g\in G:g^r=e\}$.

\begin{theorem}\label{thm:nugeneral}
For every $r\in\mathbb{Z}_{\geq1}$ with $\gcd(r,N)=1$,
\[
  \nu_r(x)=\frac k2\Bigl(1-\Jac{\mathfrak{d}}{r}\Bigr)\in\{0,k\}.
\]
\end{theorem}

\begin{proof}
First, $\exp(G)$ divides $N$. Indeed, let $g\in G$ have order $m$. By Proposition \ref{prop:local}(4), the twists of the simple summands of $I(g)$ are the values $\chi(g)$ for $\chi\in\widehat{G}$, and these exhaust $\mu_m$. So $\mu_m\subseteq\mu_N$, i.e.\ $m\mid N$; taking $m=\exp(G)$ gives $\exp(G)\mid N$. Hence $\gcd(r,n)=1$, and $g^r=e$ only for $g=e$. Since $(\mathcal{N}_G)_\mathrm{pt}\simeq\mathrm{Vec}_G$ with trivial associator (Section \ref{subsec:localpt}), $\nu_r(g)=\delta_{g^r,e}$, and therefore
\[
  \nu_r(\mathcal{N}_G)=\sum_{X\in\mathcal{O}(\mathcal{N}_G)}\dim(X)\nu_r(X)=1+d\,\nu_r(x).
\]
Next, let $\zeta_N:=e^{2\pi i/N}$. Then $d\in\mathbb{Q}(\zeta_N)$, because $d=(\dim(Y_e')-1)/k$ by Proposition \ref{prop:local}(2) and the dimensions of the simple objects of $\mathcal{Z}(\mathcal{N}_G)$ lie in $\mathbb{Q}(\zeta_N)$ (see e.g.\ \cite[Section 2]{MR3486174}). In particular $\mathbb{Q}(\sqrt\Delta)\subseteq\mathbb{Q}(\zeta_N)$, so $\mathfrak{f}$ divides $N$. Let $\tilde r$ be an inverse of $r$ modulo $N$ and $\sigma_r$ the automorphism of $\mathbb{Q}(\zeta_N)$ with $\sigma_r(\zeta_N)=\zeta_N^{\tilde r}$. By \cite[Theorem 4.4]{MR3997136} with $a=1$, and $\nu_1(\mathcal{N}_G)=1$,
\[
  \nu_r(\mathcal{N}_G)=\frac{\dim(\mathcal{N}_G)}{\sigma_r(\dim(\mathcal{N}_G))}=\frac{\alpha}{\sigma_r(\alpha)},
\]
as $\dim(\mathcal{N}_G)=n\alpha$. If $\sigma_r$ fixes $\mathbb{Q}(d)$ this is $1$. Otherwise it is $\alpha/\alpha'=\alpha^2/(4+k^2n)=\alpha-1=1+kd$, by \eqref{eq:ddprime} and \eqref{eq:ids}. Comparing with $1+d\nu_r(x)$ gives $\nu_r(x)=0$ or $k$ in the two cases. Finally, $\sigma_r$ acts on $\mathbb{Q}(\sqrt\Delta)=\mathbb{Q}(\sqrt{\mathfrak{d}})$ by the Kronecker symbol $\Jac{\mathfrak{d}}{\tilde r}$. This equals $\Jac{\mathfrak{d}}{r}$, because $\Jac{\mathfrak{d}}{\cdot}$ is a character modulo $\mathfrak{f}$, $\mathfrak{f}\mid N$ and $r\tilde r\equiv1\pmod N$.
\end{proof}

\begin{corollary}\label{cor:fsexp}
\leavevmode
\begin{enumerate}
\item A prime divides $N$ if and only if it divides $n(k^2n+4)$.
\item Both $\exp(G)$ and $\mathfrak{f}$ divide $N$.
\item If $k\geq2$, then $N$ is even.
\item If $n$ is odd, then $k=1$ or $k$ is even.
\end{enumerate}
\end{corollary}

\begin{proof}
By the Cauchy theorem \cite[Theorem~3.9]{MR3486174}, a prime ideal of $\mathbb{Z}[\zeta_N]$ divides $N$ if and only if it divides the algebraic integer $\dim(\mathcal{N}_G)=n\alpha$. So a rational prime divides $N$ if and only if it divides the norm of $n\alpha$, which is a power of $n^2\alpha\alpha'=n^2(k^2n+4)$ by \eqref{eq:ddprime}; this proves (1). That $\exp(G)$ and $\mathfrak{f}$ divide $N$ was shown in the proof of Theorem \ref{thm:nugeneral}. If $N$ were odd, then Theorem \ref{thm:nugeneral} with $r=2$ would give $\nu_2(x)\in\{0,k\}$; but $x$ is self-dual, so $\nu_2(x)=\pm1$ \cite[Theorem 2.7]{MR3427429}, and hence $k=1$, which proves (3). Finally, if $n$ and $k$ are both odd, then $n(k^2n+4)$ is odd, so $N$ is odd by (1), and therefore $k=1$ by (3).
\end{proof}

\par The indicators also constrain the twists of the simple objects of $\mathcal{Z}(\mathcal{N}_G)$ that do not lie over any invertible object. Let $\mathcal{W}$ be the set of simple objects of $\mathcal{Z}(\mathcal{N}_G)$ that are not summands of $\bigoplus_{g\in G}I(g)$. By Proposition \ref{prop:local}(5), each $W\in\mathcal{W}$ has $F(W)=m_Wx$ with $m_W\geq1$. Put $u_g:=\ell_{1,g}\ell_{2,g}$, so that $u_e=0$ and $0\leq u_g\leq k^2/4$.

\begin{lemma}\label{lem:moments}
For every $r\geq1$, with $\tau_r:=\sum_{g\in G}q(g)^r=\tau_r(\mathcal{C}(G,q))$,
\begin{align*}
E_0&:=\sum_{W\in\mathcal{W}}m_W^2=2n+\tfrac12k^2n(n-1)+2\sum_{g\in G}u_g,\\
E_r&:=\sum_{W\in\mathcal{W}}m_W^2\,\theta_W^{\,r}=\Bigl(\nu_r(x)-\tfrac k2|G[r]|\Bigr)\sqrt\Delta-\tfrac{k^2}{2}\,\tau_r+2\sum_{g\in G}u_g\,q(g)^r .
\end{align*}
\end{lemma}

\begin{proof}
By Propositions \ref{prop:local}(1)--(4) and \ref{prop:tq}, the simple summands of $\bigoplus_gI(g)$, i.e.\ the simple objects outside $\mathcal{W}$, are pairwise non-isomorphic and of two kinds. The objects of the first kind are $Y_g$ and $Y_g'$ for $g\in G$; they have twist $q(g)$, dimensions $1+\ell_{1,g}d$ and $1+\ell_{2,g}d$, and $[F(Y_g),x]=\ell_{1,g}$, $[F(Y_g'),x]=\ell_{2,g}$. The objects of the second kind are the $\frac12n(n-1)$ objects $Y_\gamma$, one for each pair $\{\gamma,\gamma^\sharp\}$ in $\Gamma\setminus\mathcal{T}$; they have dimension $\alpha$ and $[F(Y_\gamma),x]=k$, and $Y_\gamma$ has twist $\chi(g)$ for $\gamma=(g,\chi)$. Since $\sum_{(g,\chi)\in\Gamma}\chi(g)^r=n|G[r]|$, the sum of $\theta^r$ over the objects of the second kind, i.e.\ over the pairs $\{\gamma,\gamma^\sharp\}$, is $\frac12(n|G[r]|-\tau_r)$. By \eqref{eq:NS} and Frobenius reciprocity, using $\ell_{1,g}+\ell_{2,g}=k$ and $\dim(\mathcal{N}_G)=n\alpha$,
\[
\nu_r(x)\,n\alpha=\sum_{g\in G}q(g)^r\bigl(k+(k^2-2u_g)d\bigr)+\tfrac{k\alpha}{2}\bigl(n|G[r]|-\tau_r\bigr)+d\,E_r .
\]
Divide by $d$, and use $n\alpha/d=d-d'=\sqrt\Delta$ and $k+k^2d-\frac k2\alpha=\frac{k^2}{2}d$. In the same way, $\sum_X\dim(X)^2=\dim(\mathcal{N}_G)^2$ reads
\[
n\bigl(2+2kd+k^2d^2\bigr)-2d^2\sum_{g}u_g+\tfrac12n(n-1)\alpha^2+d^2E_0=n^2\alpha^2,
\]
which gives $E_0$ after using $d^2=knd+n$.
\end{proof}

\subsection{When \texorpdfstring{$k=1$}{k=1}, or under hypothesis (H)}\label{subsec:k1}

\begin{corollary}\label{cor:k1}
    Let $k=1$.  Then $K=G$ in Lemma \ref{lem:liftable}, i.e.\ $\mathcal{N}_G$ satisfies \eqref{eq:star}: the objects $Y_g$ for $g\in G$ are invertible, and $F$ restricts to an equivalence $\mathcal{Z}(\mathcal{N}_G)_\mathrm{pt}\simeq\mathrm{Vec}_G$ of fusion categories.
\end{corollary}

\begin{proof}
    By Proposition \ref{prop:local}(2), $\ell_{1,g}+\ell_{2,g}=1$, so $\min(\ell_{1,g},\ell_{2,g})=0$ for all $g\in G$. By Lemma \ref{lem:liftable}(1), the tensor functor $F$ therefore restricts to a group isomorphism $\mathcal{O}(\mathcal{Z}(\mathcal{N}_G)_\mathrm{pt})\to G$. As $(\mathcal{N}_G)_\mathrm{pt}\simeq\mathrm{Vec}_G$, this gives the desired monoidal equivalence.
\end{proof}

For the rest of this subsection we assume \eqref{eq:star}, for instance $k=1$. Then $\mathcal{Z}(\mathcal{N}_G)_\mathrm{pt}$ is a pointed braided fusion category whose group of invertible objects is $\{Y_g:g\in G\}\cong G$. A pointed braided fusion category is determined up to braided equivalence by its group of invertible objects together with the quadratic form $Y\mapsto c_{Y,Y}\in\mathbb{C}^\times$ \cite[Proposition 2.41]{MR2609644}. Here $c_{Y_g,Y_g}=\theta_{Y_g}$ because $\dim(Y_g)=1$, so $\mathcal{Z}(\mathcal{N}_G)_\mathrm{pt}\simeq\mathcal{C}(G,q)$, where $q(g)=\theta_{Y_g}$ is the quadratic form of Proposition \ref{prop:tq}, with bicharacter $b$ given by \eqref{eq:bt}. Its symmetric center is $\mathcal{C}(\operatorname{rad}(b),q|_{\operatorname{rad}(b)})$, with $\operatorname{rad}(b)$ the radical defined before Proposition \ref{cor:cyclic}.

\begin{lemma}\label{lem:basic}
Assume \eqref{eq:star}. Sending $\chi\in\widehat{G}\setminus\{1\}$ to the element $f\neq e$ with $F(Y_{(e,\chi)})=e+f+kx$ is a bijection $\widehat{G}\setminus\{1\}\to G\setminus\{e\}$; for $f\in G\setminus\{e\}$ let $\chi_f\in\widehat{G}\setminus\{1\}$ denote the character corresponding to $f$. Then:
\begin{enumerate}
\item $\chi_{f^{-1}}=\chi_f^{-1}$ for every $f\in G\setminus\{e\}$;
\item if $\chi\neq t_g$ and $(g,\chi)^\sharp=(h,\psi)$, then $\chi=t_g\,\chi_{g^{-1}h}$;
\item $\chi_f^{\,2}=b(-,f)$ for every $f\in G\setminus\{e\}$.
\end{enumerate}
\end{lemma}

\begin{proof}
The bijection is Proposition \ref{prop:local}(3) with $g=e$, since $t_e=1$. By Proposition \ref{prop:local}(6), $Y_{(e,\chi_f)}^{*}\cong Y_{(e,\chi_f^{-1})}$, and this object has forgetful image $e+f^{-1}+kx$; hence $\chi_{f^{-1}}=\chi_f^{-1}$, which is (1). Now let $\chi\neq t_g$ with $(g,\chi)^\sharp=(h,\psi)$. Lemma \ref{lem:liftable}(3) with $u=g$, which applies to every $g$ under \eqref{eq:star}, together with $t_e=1$, gives $Y_g\otimes Y_{(e,\chi')}\cong Y_{(g,t_g\chi')}$ for every $\chi'\neq1$. Put $\chi':=t_g^{-1}\chi\neq1$ and write $\chi'=\chi_{f'}$. Then $Y_{(g,\chi)}\cong Y_g\otimes Y_{(e,\chi')}$ has forgetful image $g+gf'+kx$, so $gf'=h$, that is, $\chi=t_g\,\chi_{g^{-1}h}$, which is (2). Finally, write $(e,\chi_f)^\sharp=(f,\psi)$. By Proposition \ref{prop:tq}(1) with $g=e$ and $h=f$, we have $\psi(u)=\chi_f(u)\,t_u(f)^{-1}$ for all $u\in G$, while (2), applied to the label $(f,\psi)$ whose partner is $(e,\chi_f)$, gives $\psi=t_f\,\chi_{f^{-1}}=t_f\,\chi_f^{-1}$ by (1). Comparing the two expressions for $\psi$ and using \eqref{eq:bt}, we obtain $\chi_f(u)^2=t_f(u)\,t_u(f)=b(u,f)$ for all $u\in G$, which is (3).
\end{proof}
We now determine $\operatorname{rad}(b)$ for every finite abelian group $G$: it is the $2$-torsion subgroup $G[2]$. Nothing in Subsection \ref{subsec:localpt} or Lemma \ref{lem:basic} uses the parity of $|G|$.

\begin{theorem}\label{thm:radical}
Assume \eqref{eq:star}, for instance $k=1$, and let $G$ be any finite abelian group.  Then
\[
  \operatorname{rad}(b)=G[2].
\]
\end{theorem}

\begin{proof}
By Lemma \ref{lem:basic}, $f\mapsto\chi_f$ is a bijection from $G\setminus\{e\}$ onto the set of nontrivial characters of $G$, and for all $g\in G$ and $f\neq e$,
\[
  b(g,f)=\chi_f(g)^{2}=\chi_f(g^{2}).
\]
Hence $g\in\operatorname{rad}(b)$ if and only if $\chi(g^{2})=1$ for every $\chi\in\widehat{G}$, that is, if and only if $g^{2}=e$.
\end{proof}

In particular the symmetric center $\mathcal{C}(G[2],q|_{G[2]})$ of $\mathcal{Z}(\mathcal{N}_G)_\mathrm{pt}$ has dimension $|G[2]|$, and $q$ is nondegenerate if $|G|$ is odd.

\begin{corollary}\label{cor:oddfactor}
    Assume \eqref{eq:star}, for instance $k=1$, and that $|G|$ is odd. Then $\mathcal{Z}(\mathcal{N}_G)\simeq\mathcal{Z}(\mathcal{N}_G)_\mathrm{pt}\boxtimes\mathcal{Z}(\mathcal{N}_G)_\mathrm{ad}$ is a braided equivalence and $\mathcal{O}(\mathcal{Z}(\mathcal{N}_G)_\mathrm{pt})\cong G$; that is, $\mathcal{N}_G$ is $G$-factorizable in the sense of Definition \ref{def:factorizable}.
\end{corollary}

\begin{proof}
    By Theorem \ref{thm:radical}, $\mathcal{Z}(\mathcal{N}_G)_\mathrm{pt}\simeq\mathcal{C}(G,q)$ is nondegenerate. So the factorization follows from \cite[Theorem 4.2]{MR1990929}, and the centralizer of $\mathcal{Z}(\mathcal{N}_G)_\mathrm{pt}$ is $\mathcal{Z}(\mathcal{N}_G)_\mathrm{ad}$ by \cite[Corollary 3.27]{MR2609644}. Finally, $\mathcal{O}(\mathcal{Z}(\mathcal{N}_G)_\mathrm{pt})\cong G$ is \eqref{eq:star}.
\end{proof}

The results of Subsections \ref{subsec:localpt}, \ref{subsec:prim} and \ref{subsec:indicators}, including Proposition \ref{cor:cyclic}, hold for every $k$ without \eqref{eq:star}.


\section{Factorizable near-groups}\label{sec:factorizable}

Corollary \ref{cor:oddfactor} motivates the following definition.

\begin{definition}\label{def:factorizable}
    We say that $\mathcal{N}_G$ is a \emph{$G$-factorizable near-group} if $\mathcal{N}_G$ is a $G+k|G|$ near-group fusion category for a finite abelian group $G$ and integer $k\geq1$ such that $\mathcal{Z}(\mathcal{N}_G)\simeq\mathcal{Z}(\mathcal{N}_G)_\mathrm{pt}\boxtimes\mathcal{Z}(\mathcal{N}_G)_\mathrm{ad}$ is a braided equivalence, and $\mathcal{O}(\mathcal{Z}(\mathcal{N}_G)_\mathrm{pt})\cong G$.
\end{definition}

\par Let $\mathcal{N}_G$ be a $G$-factorizable near-group. We apply the notation and conventions of Section \ref{sec:neargroup} to $\mathcal{N}_G$; in particular we may assume that $\mathcal{N}_G$ is pseudounitary, as Galois conjugation preserves $G$-factorizability, and record some consequences of Definition \ref{def:factorizable} that are used throughout this section.

\par Since $F$ is a tensor functor, it restricts to a group homomorphism from $\mathcal{O}(\mathcal{Z}(\mathcal{N}_G)_\mathrm{pt})$ to $G$. This homomorphism is injective by Lemma \ref{lem:liftable}(1), hence bijective, as $\mathcal{O}(\mathcal{Z}(\mathcal{N}_G)_\mathrm{pt})\cong G$ by definition. Thus $K=G$ in Lemma \ref{lem:liftable}, i.e.\ $\mathcal{N}_G$ satisfies \eqref{eq:star}. We may therefore take $\ell_{1,g}=0$ for all $g$, so that $Y_g$ is the invertible object with $F(Y_g)=g$, and Lemma \ref{lem:liftable}(2)--(4) applies to every element of $G$.

\par By Proposition \ref{prop:local} (see also \cite[Theorem 2.13]{MR3427429}), the simple summands of $I(e)$ are:
\begin{itemize}
\item the tensor unit $Y_e$, of dimension $1$;
\item the object $Y_e'$ with $F(Y_e')=e+kx$, of dimension $1+kd$;
\item the objects $Y_{(e,\chi)}$ for $\chi\neq1$, with $F(Y_{(e,\chi)})=e+f+kx$ for some $f\neq e$, of dimension $2+kd$.
\end{itemize}
More generally, the simple summands of $\bigoplus_{g\in G}I(g)$ are the objects $Y_g$, $Y_g'$ and $Y_{(g,\chi)}$ for $\chi\neq t_g$, of dimensions $1$, $1+kd$ and $2+kd$ respectively. By Lemma \ref{lem:liftable}(2),(3), each of these three families is closed under tensoring with invertible objects.

\par Since $d\not\in\mathbb{Q}$, we can select $\sigma\in\mathrm{Gal}(\overline{\mathbb{Q}}/\mathbb{Q})$ such that $\left.\sigma\right|_{\mathbb{Q}(d)}$ is nontrivial; we fix such a $\sigma$ for the rest of this section, and write $X\mapsto\hat\sigma(X)$ for the corresponding Galois permutation of the simple objects of the modular category $\mathcal{Z}(\mathcal{N}_G)$ \cite[Section 2]{MR3486174}.

\par The braided equivalence in Definition \ref{def:factorizable} may be taken to be the tensor product functor. Indeed, the M\"uger center of a Deligne product is the product of the M\"uger centers and $\mathcal{Z}(\mathcal{N}_G)$ is nondegenerate, so $\mathcal{Z}(\mathcal{N}_G)_\mathrm{pt}$ is nondegenerate. Hence the tensor product functor $\mathcal{Z}(\mathcal{N}_G)_\mathrm{pt}\boxtimes\mathcal{Z}(\mathcal{N}_G)_\mathrm{ad}\to\mathcal{Z}(\mathcal{N}_G)$ is a braided equivalence by \cite[Theorem 4.2]{MR1990929} and \cite[Corollary 3.27]{MR2609644}, as in the proof of Corollary \ref{cor:oddfactor}. So every simple object of $\mathcal{Z}(\mathcal{N}_G)$ is uniquely of the form $Y_u\otimes X$ with $u\in G$ and $X\in\mathcal{O}(\mathcal{Z}(\mathcal{N}_G)_\mathrm{ad})$. In particular $Y_u\otimes X\not\cong X$ for all simple $X$ and $u\neq e$, and each orbit of simple objects under tensoring with invertible objects contains exactly one object of $\mathcal{Z}(\mathcal{N}_G)_\mathrm{ad}$.


\subsection{Simple objects}

\begin{proposition}\label{factorprop}
If $\mathcal{N}_G$ is a $G$-factorizable near-group, then $|G|$ is odd and Figure \ref{fig:simples} indexes the isomorphism classes of simple objects of $\mathcal{Z}(\mathcal{N}_G)_\mathrm{ad}$, lists their dimensions, twists, and forgetful images in $\mathcal{N}_G$.
\end{proposition}

\begin{figure}[H]
    \centering
    \begin{equation*}
        \begin{array}{|c|c|c|c|c|c|}\hline
            \text{Label} & \text{Index} & \text{Quantity} & \text{Dimension} & \text{Twist} & \text{Forgetful image}  \\\hline
            Y_e,Y_e' & \text{none} & 1,1 & 1,1+kd & 1,1 & e,e+kx \\\hline
            Z_g=Y_{(g,\,t_g\overline{b(g,-)})} & g\in G_{\pm} & \frac{1}{2}(|G|-1) & 2+kd & q(g)^{-1} & g+g^{-1}+kx  \\\hline
            Z_j & j\in\mathcal{J} & |\mathcal{J}| & m_jd & \theta_{Z_j} & m_jx  \\\hline
        \end{array}
    \end{equation*}
    \caption{Isomorphism classes of simple objects of $\mathcal{Z}(\mathcal{N}_G)_\mathrm{ad}$. Here $G_{\pm}$ is the set $G\setminus\{e\}$ modulo the relation $g\sim g^{-1}$ (note $Z_g\cong Z_{g^{-1}}$), $\mathcal{J}$ is an unknown indexing set, and the twists $\theta_{Z_j}$ are not determined here.}
    \label{fig:simples}
\end{figure}

\begin{proof}
We use the consequences of Definition \ref{def:factorizable} recorded above. We first show that $|G|$ is odd. Suppose $f\in G$ had order $2$, and let $Y:=Y_{(e,\chi)}$ with $F(Y)=e+f+kx$. Then $Y_f\otimes Y$ is simple with $F(Y_f\otimes Y)=f+e+kx$, so $Y_f\otimes Y\cong Y$ by Proposition \ref{prop:local}(3). This contradicts the uniqueness of the factorization $Y_u\otimes X$. Hence $G$ has no elements of order $2$, and $|G|$ is odd by Cauchy's theorem.

\par Since $\mathcal{Z}(\mathcal{N}_G)_\mathrm{ad}$ is the centralizer of $\mathcal{Z}(\mathcal{N}_G)_\mathrm{pt}$ \cite[Corollary 3.27]{MR2609644}, Lemma \ref{lem:liftable}(4) shows that a simple summand of $U_{(h,\psi)}$ lies in $\mathcal{Z}(\mathcal{N}_G)_\mathrm{ad}$ if and only if $t_g(h)\psi(g)=1$ for all $g\in G$. By \eqref{eq:bt}, this holds if and only if
\begin{equation}\label{eq:adlabel}
  \psi=t_h\,\overline{b(h,-)}.
\end{equation}
For $\psi=t_h$, condition \eqref{eq:adlabel} means $b(h,-)\equiv1$. The factorization makes $\mathcal{Z}(\mathcal{N}_G)_\mathrm{pt}\simeq\mathcal{C}(G,q)$ nondegenerate, so this forces $h=e$. Thus $Y_e$ and $Y_e'$ are the only objects of the form $Y_h,Y_h'$ in the adjoint subcategory. For $h\neq e$, the character \eqref{eq:adlabel} differs from $t_h$, and we put
\[
Z_h:=Y_{(h,\,t_h\overline{b(h,-)})}.
\]
These are all the objects of the form $Y_{(h,\psi)}$ in the adjoint subcategory.

\par We claim that $Z_h\cong Z_{h^{-1}}$. Write $F(Z_h)=h+h'+kx$. By Proposition \ref{prop:local}(6), $Z_h^\ast$ has forgetful image $h^{-1}+h'^{-1}+kx$, and so does the simple object $Y_{(hh')^{-1}}\otimes Z_h$. Both are summands of $I(h^{-1})$ whose forgetful images contain $h'^{-1}$, so they are isomorphic by Proposition \ref{prop:local}(3). As $\mathcal{Z}(\mathcal{N}_G)_\mathrm{ad}$ is closed under duals, $Y_{(hh')^{-1}}\otimes Z_h$ lies in $\mathcal{Z}(\mathcal{N}_G)_\mathrm{ad}$ together with $Z_h$. By the uniqueness of the factorization this forces $Y_{(hh')^{-1}}\cong Y_e$, i.e.\ $h'=h^{-1}$. Consequently the partner of the label $(h,t_h\overline{b(h,-)})$ is the unique label over $h^{-1}$ of the form \eqref{eq:adlabel}, namely its inverse $(h^{-1},t_{h^{-1}}\overline{b(h^{-1},-)})$, and hence $Z_h\cong Z_{h^{-1}}\cong Z_h^\ast$. Since $|G|$ is odd, $h\neq h^{-1}$ for $h\neq e$. We index these objects by $g\in G_{\pm}$, the set $G\setminus\{e\}$ modulo the relation $g\sim g^{-1}$; there are $\frac12(|G|-1)$ of them. By Proposition \ref{prop:local}(4) and \eqref{eq:bt}, the twist of $Z_g$ is $t_g(g)\,\overline{b(g,g)}=q(g)\,q(g)^{-2}=q(g)^{-1}$, while $Y_e=\mathbbm{1}$ and $Y_e'$ are summands of $A_G\subset I(\mathbbm{1})$, so have twist $1$ (Section \ref{subsec:conventions}).

\par Any remaining simple object $Z_j\in\mathcal{Z}(\mathcal{N}_G)_\mathrm{ad}$, for some indexing $j\in\mathcal{J}$, is not a summand of any $I(g)$. By Proposition \ref{prop:local}(5), $F(Z_j)=m_jx$ and $Z_j$ has dimension $m_jd$ for some $m_j\in\mathbb{Z}_{\geq1}$, and we are done.
\end{proof}

\begin{corollary}\label{cor:factiff}
    Assume that $\mathcal{N}_G$ satisfies \eqref{eq:star}, for instance $k=1$.  Then $\mathcal{N}_G$ is $G$-factorizable if and only if $|G|$ is odd.
\end{corollary}

\begin{proof}
    If $\mathcal{N}_G$ is $G$-factorizable, then $|G|$ is odd by Proposition \ref{factorprop}.  Conversely, if $|G|$ is odd then $\mathcal{N}_G$ is $G$-factorizable by Corollary \ref{cor:oddfactor}.
\end{proof}
    
\par We now specialize the twist identities of Section \ref{subsec:indicators} to $G$-factorizable near-groups. After Definition \ref{def:factorizable} we took $\ell_{1,g}=0$, so $u_g=\ell_{1,g}\ell_{2,g}=0$ for every $g\in G$. By the factorization, the simple objects of $\mathcal{Z}(\mathcal{N}_G)$ are the $Y_u\otimes X$ with $u\in G$ and $X\in\mathcal{O}(\mathcal{Z}(\mathcal{N}_G)_\mathrm{ad})$, with $F(Y_u\otimes X)=uF(X)$, $\theta_{Y_u\otimes X}=q(u)\,\theta_X$ and $\dim(Y_u\otimes X)=\dim(X)$. By Frobenius reciprocity, $Y_u\otimes X$ is a summand of $\bigoplus_{g\in G}I(g)$ if and only if $F(X)$ contains an invertible object. So by Figure \ref{fig:simples}, the set $\mathcal{W}$ of Lemma \ref{lem:moments} consists of the objects $Y_u\otimes Z_j$ with $u\in G$ and $j\in\mathcal{J}$, and $m_{Y_u\otimes Z_j}=m_j$. Hence, in the notation of Lemma \ref{lem:moments},
\begin{equation}\label{eq:Wfact}
  E_0=|G|\sum_{j\in\mathcal{J}}m_j^2,\qquad E_r=\tau_r\sum_{j\in\mathcal{J}}m_j^2\,\theta_{Z_j}^{\,r}\quad(r\geq1).
\end{equation}

\begin{proposition}\label{prop:sumbj}
    If $\mathcal{N}_G$ is a $G$-factorizable near-group, then
\begin{equation}
    \sum_{j\in\mathcal{J}}m_j^2=\frac{1}{2}(4+(|G|-1)k^2).
\end{equation}
\end{proposition}

\begin{proof}
By \eqref{eq:Wfact} and Lemma \ref{lem:moments} with $u_g=0$, $|G|\sum_jm_j^2=2|G|+\frac12k^2|G|(|G|-1)$; divide by $|G|$.
\end{proof}

The exotic multiplicities $m_j$ are further constrained by formal codegrees.

\begin{corollary}\label{cor:codeg}
If $\mathcal{N}_G$ is a $G$-factorizable near-group, then $m_j^2$ divides $k^2|G|+4$ for every $j\in\mathcal{J}$. Consequently:
\begin{enumerate}
\item $\gcd(m_j,k)\in\{1,2\}$ for every $j\in\mathcal{J}$; in particular, if $k\geq3$ then $F(Z_j)\neq kx$ for every $j\in\mathcal{J}$;
\item if $k=1$ and $|G|+4$ is squarefree, then $F(Z_j)=x$ for every $j\in\mathcal{J}$, so $|\mathcal{J}|=\frac12(|G|+3)$ and $\mathcal{Z}(\mathcal{N}_G)_\mathrm{ad}$ has rank $|G|+3$.
\end{enumerate}
\end{corollary}

\begin{proof}
The category $\mathcal{Z}(\mathcal{N}_G)_\mathrm{ad}$ is modular, being a Deligne factor of the modular category $\mathcal{Z}(\mathcal{N}_G)$, and its global dimension is $\dim(\mathcal{Z}(\mathcal{N}_G))/\dim(\mathcal{Z}(\mathcal{N}_G)_\mathrm{pt})=|G|\alpha^2$. For a simple object $X$ of a modular category $\mathcal{M}$, the assignment $Y\mapsto S_{Y,X}/S_{\mathbbm{1},X}$ is a character of the Grothendieck ring of $\mathcal{M}$ by the Verlinde formula, and its formal codegree is $\dim(\mathcal{M})/\dim(X)^2$. Formal codegrees are algebraic integers \cite[Theorem 1.2]{MR2576705}. For $X=Z_j$, by \eqref{eq:ids} this formal codegree is
\[
\frac{|G|\alpha^2}{m_j^2d^2}=\frac{\alpha^2}{m_j^2(\alpha-1)}=\frac{k^2|G|+4}{m_j^2},
\]
a rational algebraic integer, hence an integer.

For (1), let $p$ be a prime dividing both $m_j$ and $k$. Then $p^2$ divides $k^2|G|+4$ and $k^2|G|$, hence divides $4$, so $p=2$. If $4$ divided both $m_j$ and $k$, then $16$ would divide $k^2|G|+4$ and $k^2|G|$, which is impossible. So $\gcd(m_j,k)\in\{1,2\}$, and $m_j=k$ forces $k\leq2$. For (2), $m_j^2$ divides the squarefree integer $|G|+4$, so $m_j=1$ for every $j\in\mathcal{J}$; then $|\mathcal{J}|=\sum_jm_j^2=\frac12(|G|+3)$ by Proposition \ref{prop:sumbj}, and by Figure \ref{fig:simples} there are $2+\frac12(|G|-1)$ simple objects besides the $Z_j$.
\end{proof}

\begin{note}\label{rem:codeg}
Since $|G|$ is odd for a $G$-factorizable near-group, $k$ is $1$ or even by Corollary \ref{cor:fsexp}(4), and for even $k$ the integer $k^2|G|+4$ is divisible by $4$. So Corollary \ref{cor:codeg}(2) covers every case in which $k^2|G|+4$ is squarefree. For $k=1$ Corollary \ref{cor:codeg} gives more. As $|G|+4$ is odd, every $m_j$ is odd. As $m_j^2\leq\sum_{i\in\mathcal{J}}m_i^2=\frac12(|G|+3)<|G|+4$ by Proposition \ref{prop:sumbj}, $m_j^2$ is a proper divisor of the odd integer $|G|+4$, so $m_j^2\leq\frac13(|G|+4)$. In particular, if no square greater than $1$ properly divides $|G|+4$ (for instance, if $|G|+4$ is squarefree or the square of a prime), then $m_j=1$ for every $j\in\mathcal{J}$. Finally, odd squares are $\equiv1\pmod8$, so $|\mathcal{J}|\equiv\frac12(|G|+3)\pmod 8$. This is consistent with the modular data for doubles of near-groups of type $G+|G|$ \cite{MR3167494}.
\end{note}

\subsection{Fusion rules of \texorpdfstring{$\mathcal{Z}(\mathcal{N}_G)_\mathrm{ad}$}{Z(N\_G)-ad}}

\begin{proposition}\label{fus}
Assume $\mathcal{N}_G$ is a $G$-factorizable near-group, and let $g\in G_{\pm}$ and $i\in\mathcal{J}$. Then the first, third and fourth identities below hold, and the second holds for every $h\in G_{\pm}$ with $h\neq g$ (the classes of $gh$ and $gh^{-1}$ in $G_{\pm}$ do not depend on the chosen representatives of $g$ and $h$):
\begin{align}
&Z_g\otimes Y_e'=k^2Y_e'+(k^2+1)Z_g+k^2\sum_{\substack{f\in G_{\pm}\\ f\neq g}}Z_f+k\sum_{j\in\mathcal{J}}m_jZ_j \\
&Z_g\otimes Z_h=k^2Y_e'+(k^2+1)\left(Z_{gh}+Z_{gh^{-1}}\right)+k^2\sum_{\substack{ f\in G_{\pm} \\ f\neq gh,gh^{-1}}}Z_f+k\sum_{j\in\mathcal{J}}m_jZ_j \\
&Z_g\otimes Z_g=Y_e+(k^2+1)(Y_e'+Z_{g^2})+k^2\sum_{\substack{ f\in G_{\pm} \\ f\neq g^2}}Z_f+k\sum_{j\in\mathcal{J}}m_jZ_j \\
&Z_g\otimes Z_i=m_ikY_e'+m_ik\sum_{f\in G_{\pm}}Z_f+m_i\sum_{j\in\mathcal{J}}m_jZ_j.
\end{align}
\end{proposition}

\begin{proof}
First note that for $g\in G_{\pm}$, the multiplicity of $Z_g$ in any tensor product is determined by the multiplicity of $g$ since this is the only simple summand of $I(g)$ lying in the adjoint.  This determines most fusion rules, e.g. $[Z_j,Z_g\otimes Y_e']=[Y_e'\otimes Z_j,Z_g]=m_jk$ since
\begin{equation}
F(Y_e'\otimes Z_j)=(e+kx)m_jx=m_jk\sum_{g\in G}g+(1+k^2|G|)m_jx.    
\end{equation}
Similarly we have $[Y_e',Z_g\otimes Y_e']=[Z_g,Y_e'\otimes Y_e']=[g,F(Y_e')^2]=k^2$, as $F(Y_e')^2=(e+kx)^2$.  Now,
\begin{equation}
F(Z_g\otimes Y_e')=(g+g^{-1}+kx)(e+kx)=g+g^{-1}+k^2\sum_{h\in G}h+k(3+k^2|G|)x.
\end{equation}
Therefore $[Z_h,Z_g\otimes Y_e']=k^2$ for all $h\neq g$ in $G_{\pm}$ and $[Z_g,Z_g\otimes Y_e']=k^2+1$.  The simple objects $Z_g$ are self-dual for all $g\in G_{\pm}$ (see the proof of Proposition \ref{factorprop}). So we compute for all $h\neq g$ in $G_{\pm}$,
\begin{align}
F(Z_g\otimes Z_h)&=(g+g^{-1}+kx)(h+h^{-1}+kx) \\
&=gh+g^{-1}h^{-1}+gh^{-1}+g^{-1}h+4kx+ k^2\left(k|G|x+\sum_{f\in G}f\right) \\
&=gh+ g^{-1}h^{-1}+ gh^{-1}+ g^{-1}h+k^2\sum_{f\in G}f+k(4+k^2|G|)x
\end{align}
Thus for all $f$ distinct from $gh$ and $gh^{-1}$ in $G_{\pm}$, $[Z_f,Z_g\otimes Z_h]=k^2$ while $[Z_f,Z_g\otimes Z_h]=k^2+1$ if $f=gh$ or $f=gh^{-1}$ in $G_{\pm}$.  Also $[Y_e,Z_g\otimes Z_h]=0$ and $[Y_e',Z_g\otimes Z_h]=k^2$.  We compute for all $j$,
\begin{equation}
F(Z_j\otimes Z_h)=m_jx(h+ h^{-1}+kx)=m_jk\sum_{g\in G}g+(2+k^2|G|)m_jx.
\end{equation}
Thus for all $h\neq g$ in $G_{\pm}$, $[Z_j,Z_g\otimes Z_h]=[Z_g,Z_j\otimes Z_h]=m_jk$ giving a complete decomposition of $Z_g\otimes Z_h$ into simple objects.  Now for any $i,j$,
\begin{equation}
F(Z_i\otimes Z_j)=m_im_jx^2=m_im_j\sum_{g\in G}g+k|G|m_im_jx
\end{equation}
thus $[Z_j,Z_g\otimes Z_i]=[Z_g,Z_j\otimes Z_{i^\ast}]=m_jm_{i^\ast}=m_im_j$, using $m_{i^\ast}=m_i$.
\end{proof}

\begin{proposition}\label{zfusion}
    Assume $\mathcal{N}_G$ is a $G$-factorizable near-group.  Then for all $i,j\in\mathcal{J}$, $Z_i\otimes Z_j$ is
    \begin{equation}
\left\{\begin{array}{ccc}
Y_e+(m_im_j-1)Y_e'+m_im_j\displaystyle{\sum_{g\in G_{\pm}}}Z_g+\displaystyle{\sum_{\ell\in\mathcal{J}_{i,j}}}m_{i,j,\ell}Z_\ell & : & i^\ast=j\\
m_im_jY_e'+m_im_j\displaystyle{\sum_{g\in G_{\pm}}}Z_g+\displaystyle{\sum_{\ell\in\mathcal{J}_{i,j}}}m_{i,j,\ell}Z_\ell &:  & i^\ast\neq j \\
\end{array}
\right.
\end{equation}
for some unknown indexing sets $\mathcal{J}_{i,j}\subset\mathcal{J}$ and multiplicities $m_{i,j,\ell}\in\mathbb{Z}_{\geq1}$. These satisfy
\begin{equation}\label{eq:mijl}
\sum_{\ell\in\mathcal{J}_{i,j}}m_{i,j,\ell}\,m_\ell=k\Bigl(\tfrac12(|G|-1)\,m_im_j+\delta_{i^\ast,j}\Bigr).
\end{equation}
\end{proposition}

\begin{proof}
For any $i,j\in\mathcal{J}$,
\begin{equation}
F(Z_i\otimes Z_j)=m_im_j\sum_{g\in G}g+km_im_j|G|x.
\end{equation}
The multiplicities of $Y_e$, $Y_e'$, and $Z_g$ for $g\in G_{\pm}$ are determined by the forgetful functor, the multiplicity of $Y_e$ being $1$ if $j=i^\ast$ and $0$ otherwise; thus for some indexing sets $\mathcal{J}_{i,j}\subset\mathcal{J}$ and multiplicities $m_{i,j,\ell}$ for $\ell\in \mathcal{J}_{i,j}$ we have our result. Finally, the multiplicity of $x$ in $F(Z_i\otimes Z_j)$ is $km_im_j|G|$, while by Figure \ref{fig:simples} the decomposition above contributes $k(m_im_j-\delta_{i^\ast,j})+\tfrac k2(|G|-1)m_im_j+\sum_\ell m_{i,j,\ell}m_\ell$; comparing gives \eqref{eq:mijl}.
\end{proof}

\begin{proposition}\label{prop:yeprime}
    Assume $\mathcal{N}_G$ is a $G$-factorizable near-group.  Then
    \begin{align}
        Y_e'\otimes Y_e'&=Y_e+k^2Y_e'+k^2\sum_{g\in G_{\pm}}Z_g+k\sum_{j\in\mathcal{J}}m_jZ_j \\
        Y_e'\otimes Z_j&=km_jY_e'+km_j\sum_{g\in G_{\pm}}Z_g+m_j\sum_{\substack{\ell\in\mathcal{J}\\ \ell\neq j}}m_\ell Z_\ell+(m_j^2-1)Z_j.
    \end{align}
\end{proposition}

\begin{proof}
    We compute
\begin{equation}\label{yeprime}
    F(Y_e'\otimes Z_j)=(e+kx)(m_jx)=m_jx+km_jx^2=km_j\sum_{g\in G}g+m_j(1+k^2|G|)x
\end{equation}
hence $[Z_g,Y_e'\otimes Z_j]=km_j$ for all $g\in G_{\pm}$ and we already knew that $[Y_e',Y_e'\otimes Z_j]=km_j$ for all $j\in\mathcal{J}$.  Similarly,
\begin{equation}
    F(Y_e'\otimes Y_e')=(e+kx)(e+kx)=(k^2+1)e+k^2\sum_{g\neq e}g+k(2+|G|k^2)x
\end{equation}
This determines the coefficients of $Z_g$ for $g\in G_{\pm}$, $Y_e$, and $Y_e'$, since $Y_e'\cong {Y_e'}^\ast$.  We have $[e,F(Y_e'\otimes Z_j)]=m_jk$ hence $[Z_j,Y_e'\otimes Y_e']=m_jk$ for all $j\in\mathcal{J}$. The $Z_\ell$ coefficients follow similarly by Frobenius reciprocity and Proposition \ref{zfusion}.
\end{proof}


\subsection{Modular data of \texorpdfstring{$\mathcal{Z}(\mathcal{N}_G)_\mathrm{ad}$}{Z(N\_G)-ad}}\label{subsec:modulardata}

Assume for the remainder of this section that $\mathcal{N}_G$ is a $G$-factorizable near-group category. Recall $\alpha=2+kd$, and put $\lambda:=\sqrt{4+k^2|G|}$. By \eqref{eq:ids}, $\alpha^2|G|=\lambda^2d^2$, so $\alpha\sqrt{|G|}=d\lambda$.  We abbreviate the twists $\eta_j:=\theta_{Z_j}$ for $j\in\mathcal{J}$ in this section; the letter $\xi$ is reserved for the central charge introduced below. Since there is a unique simple object of $\mathcal{Z}(\mathcal{N}_G)_\mathrm{ad}$ corresponding to $g\in G_{\pm}$ and $j\in\mathcal{J}$, we will abbreviate $S_{g,j}:=S_{Z_g,Z_j}$, etc., in this section and use the indices $e$ and $e'$ for $Y_e$ and $Y_e'$. The $S$-matrix is normalized as in Section \ref{subsec:conventions}. We have $Y_e'=\hat\sigma(Y_e)$, with $\sigma$ as fixed after Definition \ref{def:factorizable}. Indeed, for every invertible $Y_u$ we have $S_{Y_u,\hat\sigma(Y_e)}/S_{Y_e,\hat\sigma(Y_e)}=\sigma(S_{Y_u,Y_e}/S_{Y_e,Y_e})=1$ \cite[Section 2]{MR3486174}, so $\hat\sigma(Y_e)$ centralizes $\mathcal{Z}(\mathcal{N}_G)_\mathrm{pt}$ and lies in $\mathcal{Z}(\mathcal{N}_G)_\mathrm{ad}$; moreover $\dim(\hat\sigma(Y_e))^2=\dim(\mathcal{Z}(\mathcal{N}_G))/\sigma(\dim(\mathcal{Z}(\mathcal{N}_G)))=(\alpha/\alpha')^2=(\alpha-1)^2$ by \eqref{eq:ids} and \eqref{eq:ddprime}, and $Y_e'$ is the only simple object of $\mathcal{Z}(\mathcal{N}_G)_\mathrm{ad}$ of dimension $\pm(\alpha-1)$ (Figure \ref{fig:simples}). Hence for all simple $X$ in $\mathcal{O}(\mathcal{Z}(\mathcal{N}_G)_\mathrm{ad})$,
\begin{equation}
\sigma\left(\dfrac{S_{e,X}}{S_{e,e}}\right)=\dfrac{S_{e',X}}{S_{e',e}}
\end{equation}
and therefore
\begin{equation}
S_{e',X}=(\alpha-1)\sigma\left(\dim(X)\right)=\left\{\begin{array}{ccl}
\alpha-1 & : & \dim(X)=1 \\
1 & : & \dim(X)=\alpha-1 \\
\alpha & : & \dim(X)=\alpha \\
-m_jd & : & \dim(X)=m_jd
\end{array}\right..
\end{equation}

By the factorization, $\mathcal{Z}(\mathcal{N}_G)_\mathrm{pt}\simeq\mathcal{C}(G,q)$ is nondegenerate, with $q$ as in Proposition \ref{prop:tq}; recall that $q(g^{-1})=q(g)$, and that $\theta_{Z_g}=q(g)^{-1}$ (Figure \ref{fig:simples}). Hence $\sum_{g\in G_{\pm}}\theta_{Z_g}=\frac{1}{2}(\xi^{-1}\sqrt{|G|}-1)$ where $\xi=\frac{1}{\sqrt{|G|}}\sum_{g\in G}q(g)$ is the central charge of $\mathcal{Z}(\mathcal{N}_G)_\mathrm{pt}$.

\par The entries $S_{g,j}$ and $S_{g,h}$ for $g,h\in G_{\pm}$ and $j\in\mathcal{J}$ follow directly from the branching rule (Theorem \ref{thm:branch}). For $g\in G\setminus\{e\}$ write $\psi_g:=t_g\,\overline{b(g,-)}$, so that $Z_g=Y_{(g,\psi_g)}$ (Figure \ref{fig:simples}); by the proof of Proposition \ref{factorprop}, the two labels of $Z_g$ are $(g,\psi_g)$ and $(g,\psi_g)^\sharp=(g^{-1},\psi_g^{-1})$. Since $Z_j$ is not a summand of any $I(f)$, Corollary \ref{cor:branch}(1) gives
\begin{equation}\label{eq:Sgj}
S_{g,j}=0\qquad\text{for all }g\in G_{\pm},\ j\in\mathcal{J}.
\end{equation}
For $g,h\in G_{\pm}$, apply Theorem \ref{thm:branch} with $X=Z_g$ and $M=M_{(h,\psi_h)}$, so that $U(M)=Z_h$. The only simple local modules $M'$ with $[Z_g,U(M')]\neq0$ are $M_{(g,\psi_g)}$ and $M_{(g^{-1},\psi_g^{-1})}$, each with $[Z_g,U(M')]=1$ (Proposition \ref{prop:local}(3)), so by \eqref{eq:S0},
\[
S_{g,h}=\alpha\Bigl(S^0_{(g,\psi_g),(h,\psi_h)}+S^0_{(g^{-1},\psi_g^{-1}),(h,\psi_h)}\Bigr)=\alpha\bigl(\overline{\psi_g(h)\psi_h(g)}+\psi_g(h)\psi_h(g)\bigr).
\]
By \eqref{eq:bt}, $\psi_g(h)\psi_h(g)=t_g(h)\,t_h(g)\,b(g,h)^{-2}=b(g,h)^{-1}$. Hence
\begin{equation}\label{eq:Sgh}
S_{g,h}=\bigl(b(g,h)+b(g,h)^{-1}\bigr)\alpha\qquad\text{for all }g,h\in G_{\pm},
\end{equation}
including $g=h$. This does not depend on the choice of representatives of $g,h\in G_{\pm}$ in $G$, as $b(g^{-1},h)=b(g,h)^{-1}$.

Recall the unknown fusion coefficients $m_{i,j,\ell}$ from Proposition \ref{zfusion}. For $i,j\in\mathcal{J}$, the object $Z_{i^\ast}\otimes Z_j$ contains $Y_e$ with multiplicity $\delta_{i,j}$ and $Y_e'$ with multiplicity $m_im_j-\delta_{i,j}$. As $\eta_{i^\ast}=\eta_i$ and $m_{i^\ast}=m_i$, the balancing equation \eqref{eq:balancing} gives
\[
\eta_i\eta_jS_{i,j}=\delta_{i,j}+(m_im_j-\delta_{i,j})(\alpha-1)+\frac12m_im_j\alpha\bigl(\xi^{-1}\sqrt{|G|}-1\bigr)+d\sum_{\ell\in\mathcal{J}_{i^\ast,j}}m_{i^\ast,j,\ell}m_\ell\eta_\ell .
\]
Using $\alpha-2=kd$ and $\alpha\sqrt{|G|}=d\lambda$, this becomes
\begin{equation}\label{sij1}
S_{i,j}=\eta_i^{-1}\eta_j^{-1}d\left(\frac{1}{2}m_im_j\left(k+\xi^{-1}\lambda\right)-k\,\delta_{i,j}+\sum_{\ell\in\mathcal{J}_{i^\ast,j}}m_{i^\ast,j,\ell}m_\ell\eta_\ell\right).
\end{equation}
The $S$-matrix of $\mathcal{Z}(\mathcal{N}_G)_\mathrm{ad}$ is then
\begin{equation}
S = \left[\begin{array}{cc|ccc|ccc}
1 & \alpha-1 & \cdots & \alpha & \cdots & \cdots & m_jd & \cdots \\ 
\alpha-1 & 1 & \cdots & \alpha & \cdots & \cdots & -m_jd & \cdots \\\hline
\vdots & \vdots & & & & & \\
\alpha & \alpha & & (b(g,h)+b(g,h)^{-1})\alpha &  &  & 0 &  \\
\vdots & \vdots & & & & \\\hline
\vdots & \vdots &  &  &  & &    \\
m_jd & -m_jd &  &  0& &  &  S_{i,j}     \\
\vdots & \vdots &  & &  &
\end{array}\right]
\end{equation}
where the block indexed by $G_{\pm}$ is given by \eqref{eq:Sgh}, and the entries $S_{i,j}$ for $i,j\in\mathcal{J}$ are given by \eqref{sij1}.


\subsection{Gauss sums and the twists of the \texorpdfstring{$Z_j$}{Zj}}\label{subsec:twistmoments}

By \eqref{eq:Wfact}, $\sum_jm_j^2\eta_j^{\,r}=E_r/\tau_r$, where $\tau_r=\tau_r(\mathcal{C}(G,q))=\sum_{g\in G}q(g)^r$. We show that $\tau_r\neq0$ (Lemma \ref{lem:taun}) and compute it when $\gcd(r,|G|)=1$ (Proposition \ref{prop:jacobi}).

\begin{lemma}\label{lem:taun}
Let $G$ be a finite abelian group of odd order, $q$ a nondegenerate quadratic
form on $G$, and $r\in\mathbb{Z}_{\geq1}$. Then
\[
  \bigl|\tau_r\bigl(\mathcal{C}(G,q)\bigr)\bigr|^2=|G|\,|G[r]|.
\]
In particular $\tau_r(\mathcal{C}(G,q))\neq0$.
\end{lemma}

\begin{proof}
As $|G|$ is odd, $q(g)=\beta(g,g)$ for a bicharacter $\beta$ \cite[Lemma A.1]{MR2609644}, so $\mathcal{C}(G,q)$ is $\mathrm{Vec}_G$ with trivial associativity constraint and braiding $\beta$ \cite[Proposition 2.41]{MR2609644}; it is modular as $q$ is nondegenerate. Hence $\nu_r(g)=\delta_{g^r,e}$ for $g\in G$, and $\nu_r(\mathcal{C}(G,q))=|G[r]|$. So $|\tau_r(\mathcal{C}(G,q))|^2=\dim(\mathcal{C}(G,q))\,\nu_r(\mathcal{C}(G,q))=|G|\,|G[r]|$ by \cite[Proposition 3.3]{MR3997136}.
\end{proof}

\begin{proposition}\label{prop:jacobi}
Let $G$ be a finite abelian group of odd order, $q$ a nondegenerate quadratic
form on $G$, and $r\in\mathbb{Z}_{\geq1}$ with $\gcd(r,|G|)=1$.  Then
\begin{equation}\label{eq:jacobi}
  \tau_r\bigl(\mathcal{C}(G,q)\bigr)=\Jac{r}{|G|}\,\tau_1\bigl(\mathcal{C}(G,q)\bigr)
  =\Jac{r}{|G|}\,\xi\sqrt{|G|},
\end{equation}
where $\Jac{r}{|G|}$ is the Jacobi symbol and $\xi=\tau_1(\mathcal{C}(G,q))/\sqrt{|G|}$ is the central charge.
\end{proposition}

\begin{proof}
Since $|G|$ is odd, $(G,q)$ is an orthogonal direct sum of cyclic metric groups $(C_{n_i},q_i)$ with $q_i(j)=\exp(2\pi i\,a_ij^2/n_i)$ and $\gcd(a_i,n_i)=1$ \cite[Section 5, Theorem 4, and Section 6]{MR156890}, and Gauss sums are multiplicative over orthogonal direct sums, as $q(g_1g_2)=q_1(g_1)q_2(g_2)$. For odd $n$ and $\gcd(c,n)=1$, the classical evaluation of quadratic Gauss sums gives $\sum_{j\in C_n}\exp(2\pi i\,cj^2/n)=\Jac{c}{n}\sum_{j\in C_n}\exp(2\pi i\,j^2/n)$ \cite[Chapter 1]{MR1625181}. Applying this with $c=a_ir$ and $c=a_i$ gives $\tau_r(\mathcal{C}(C_{n_i},q_i))=\Jac{r}{n_i}\tau_1(\mathcal{C}(C_{n_i},q_i))$, and $\prod_i\Jac{r}{n_i}=\Jac{r}{|G|}$. For $|G|$ prime this is \cite[Example 3.15]{MR3997136}.
\end{proof}

For $\gcd(r,|G|)=1$ we abbreviate
$\varepsilon_r:=\Jac{r}{|G|}=\pm1$, so that $\varepsilon_1=1$ and, by the second supplementary law of quadratic reciprocity for the Jacobi symbol, $\varepsilon_2=(-1)^{(|G|^2-1)/8}$.

\begin{theorem}\label{thm:master}
Let $\mathcal{N}_G$ be a $G$-factorizable near-group of type $G+k|G|$ and let $r\in\mathbb{Z}_{\geq1}$, with $\tau_r=\tau_r(\mathcal{C}(G,q))$ as above. Then
\begin{equation}\label{eq:mastergen}
  \sum_{j\in\mathcal{J}}m_j^2\,\eta_j^{\,r}
  \;=\;-\frac{k^2}{2}\;+\;\Bigl(\frac{\nu_r(x)}{|G[r]|}-\frac{k}{2}\Bigr)\,\overline{\tau_r}\,\frac{\lambda}{\sqrt{|G|}} .
\end{equation}
If $\gcd(r,|G|)=1$, this reads
\begin{equation}\label{eq:master}
  \sum_{j\in\mathcal{J}}m_j^2\,\eta_j^{\,r}
  \;=\;-\frac{k^2}{2}\;+\;\varepsilon_r\Bigl(\nu_r(x)-\frac{k}{2}\Bigr)\xi^{-1}\lambda .
\end{equation}
\end{theorem}

\begin{proof}
By \eqref{eq:Wfact} and Lemma \ref{lem:moments} with $u_g=0$,
\[
  \tau_r\sum_{j\in\mathcal{J}}m_j^2\eta_j^{\,r}=\Bigl(\nu_r(x)-\frac k2|G[r]|\Bigr)\sqrt\Delta-\frac{k^2}{2}\,\tau_r .
\]
Divide by $\tau_r$, which is nonzero with $1/\tau_r=\overline{\tau_r}/(|G|\,|G[r]|)$ by Lemma \ref{lem:taun}, and use $\sqrt\Delta=\sqrt{|G|}\,\lambda$; this gives \eqref{eq:mastergen}. If $\gcd(r,|G|)=1$, then $G[r]=\{e\}$, and $\overline{\tau_r}=\varepsilon_r\,\xi^{-1}\sqrt{|G|}$ by Proposition \ref{prop:jacobi}, as $\varepsilon_r$ is real and $\overline{\xi}=\xi^{-1}$; substituting gives \eqref{eq:master}.
\end{proof}

\begin{corollary}\label{cor:n12}
Taking $r=1$ in \eqref{eq:master}, we have $\varepsilon_1=1$ and $\nu_1(x)=0$,
so
\[
  \sum_{j\in\mathcal{J}}m_j^2\,\eta_j=-\frac{k}{2}\bigl(k+\xi^{-1}\lambda\bigr).
\]
Taking $r=2$,
\[
  \sum_{j\in\mathcal{J}}m_j^2\,\eta_j^{\,2}
  =-\frac{k^2}{2}+\frac{\varepsilon_2}{2}\bigl(2\nu_2(x)-k\bigr)\xi^{-1}\lambda .
\]
\end{corollary}

For $r$ prime to $N$, Theorem \ref{thm:nugeneral} determines the remaining unknown $\nu_r(x)$.

\begin{corollary}\label{prop:nun}
Let $N$ and $\mathfrak{d}$ be as in Section \ref{subsec:indicators}, and let $r\in\mathbb{Z}_{\geq1}$ with $\gcd(r,N)=1$. Then
\begin{equation}\label{eq:nun}
  \sum_{j\in\mathcal{J}}m_j^2\,\eta_j^{\,r}
  =-\frac{k}{2}\Bigl(k+\Jac{\mathfrak{d}}{r}\,\varepsilon_r\,\xi^{-1}\lambda\Bigr).
\end{equation}
If $k\geq2$, then $N$ is even, so $r=2$ is not covered, and $\nu_2(x)=\pm1$ remains unknown in Corollary \ref{cor:n12}.
\end{corollary}

\begin{proof}
By Corollary \ref{cor:fsexp}(2) and (3), $\exp(G)$ divides $N$, so $\gcd(r,|G|)=1$, and $N$ is even if $k\geq2$. Substituting $\nu_r(x)-\frac k2=-\frac k2\Jac{\mathfrak{d}}{r}$ from Theorem \ref{thm:nugeneral} into \eqref{eq:master} gives \eqref{eq:nun}.
\end{proof}

\begin{note}\label{rem:nunGalois}
For $\gcd(r,N)=1$, the identity \eqref{eq:nun} is the image of the case $r=1$ under the automorphism $\sigma_r$ of $\mathbb{Q}(\zeta_N)$ with $\sigma_r(\zeta_N)=\zeta_N^r$. Indeed $\eta_j\in\mu_N$ and $m_j\in\mathbb{Z}$, while $\xi^{-1}\lambda=\overline{\tau_1}\sqrt\Delta/|G|$, with $\sigma_r(\overline{\tau_1})=\overline{\tau_r}=\varepsilon_r\overline{\tau_1}$ by Proposition \ref{prop:jacobi} and $\sigma_r(\sqrt\Delta)=\Jac{\mathfrak{d}}{r}\sqrt\Delta$. So Corollary \ref{prop:nun} is a consistency check, in effect of Theorem \ref{thm:nugeneral} in the factorizable case. The identities \eqref{eq:mastergen} that carry information beyond $r=1$ are those with $\gcd(r,N)>1$, such as $r=2$ when $k\geq2$.
\end{note}

\begin{example}\label{ex:C3k2}
$G$-factorizable near-groups occur beyond $k=1$. By \cite[Theorem 3.6 and its proof]{yu2023center}, there is a near-group fusion category $\mathcal{N}_G$ of type $C_3+6$ with
\[
\mathcal{Z}(\mathcal{N}_G)\simeq\mathcal{C}(C_3,q)\boxtimes\mathcal{C}(\mathfrak{sl}_3,9)^0_{C_3},\qquad q(j)=e^{2\pi ij^2/3},
\]
where the second factor is the category of local modules over the simple current algebra of $\mathcal{C}(\mathfrak{sl}_3,9)$. That factor has no nontrivial invertible objects \cite[Lemma 3.1]{yu2023center}, so $\mathcal{N}_G$ satisfies \eqref{eq:star} and is $C_3$-factorizable with $k=2$, and the second factor is $\mathcal{Z}(\mathcal{N}_G)_\mathrm{ad}$. Comparing \cite[Lemmas 3.1 and 3.4]{yu2023center} with Figure \ref{fig:simples}, the set $\mathcal{J}$ has six elements, $m_j=1$ for all $j\in\mathcal{J}$, and the twists $\eta_j$ are $i,i,i,i,-1,-1$. Here $\lambda=4$, $\xi=i$ and $\varepsilon_2=-1$. Thus $\sum_j\eta_j=-2+4i$ and $\sum_j\eta_j^2=-2$, in agreement with Corollary \ref{cor:n12}, where the second identity gives $\nu_2(x)=1$. Moreover $N=12=\mathfrak{f}$.
\end{example}


\section{De-equivariantization for cyclic \texorpdfstring{$G$}{G}}\label{sec:deeq}
 
Throughout this section $G=C_n$ is cyclic of order
$n$, written additively via an identification $G=\mathbb{Z}/n\mathbb{Z}$, $k\ge1$ is arbitrary, and we assume
hypothesis \eqref{eq:star} of Definition \ref{hyp:H},
writing $Y_g$ for the invertible object with $F(Y_g)=g$.  We do not assume the
parity of $n$, nor that $\mathcal{N}_G$ is $G$-factorizable.  When $k=1$,
\eqref{eq:star} holds by Corollary \ref{cor:k1}.

\par De-equivariantizations of near-group categories by subgroups of $G$ were studied systematically by Izumi in the unitary setting, principally for $k=1$ \cite[Section 12]{MR3635673}: for a subgroup $H\subset H^\perp$ the fusion rules of the de-equivariantization are determined in \cite[Theorems 12.3 and 12.5, Corollary 12.7]{MR3635673}, and the resulting relation between Drinfeld centers is discussed in \cite[Remark 12.10]{MR3635673}.  In particular, \cite[Example 12.14]{MR3635673} obtains categories with the fusion rules of the Haagerup category from the near-group categories for $C_9$ (compare Theorem \ref{thm:haagerup} with $m=3$).  The arguments below are purely fusion-categorical: they do not assume unitarity and apply to every $k\in\mathbb{Z}_{\geq1}$.

We use the notation of Proposition \ref{cor:cyclic}: $\kappa=q(1)$ is a primitive $n$-th root of unity, and by \eqref{eq:qbcyclic}, $t_i(j)=\kappa^{ij}$, $q(j)=\kappa^{j^2}$ and $b(i,j)=\kappa^{2ij}$ for $i,j\in\mathbb{Z}/n\mathbb{Z}$.
We write $\alpha=2+kd$ as in Section \ref{sec:neargroup}, and use the identities \eqref{eq:ids}.

 
\subsection{The maximal Tannakian subcategory}\label{subsec:maxtann}
 
For each divisor $r\mid n$ let $H_r\subset C_n$ denote the unique subgroup of order
$r$, namely $H_r=\tfrac nr C_n$.  Every fusion subcategory of $\mathcal{C}(G,q)$
is of the form $\mathcal{C}(H,q|_H)$ for a subgroup $H\subset G$ \cite[Section 2.11]{MR2609644},
so the fusion subcategories of $\mathcal{Z}(\mathcal{N}_G)_{\mathrm{pt}}$ are
exactly the $\mathcal{C}(H_r,q|_{H_r})$ for $r\mid n$. As $(H,q|_H)$ is the pre-metric group of $\mathcal{C}(H,q|_H)$ (Section \ref{subsec:k1}), $\mathcal{C}(H,q|_H)$ is Tannakian if and only if $q|_H\equiv1$ \cite[Example 2.48]{MR2609644}, and then $\mathcal{C}(H,q|_H)\simeq\mathrm{Rep}(\widehat{H})$ \cite[Example 2.45]{MR2609644}.
 
\begin{proposition}\label{prop:maxtann}
Let $r\mid n$.  Then $\mathcal{C}(H_r,q|_{H_r})$ is Tannakian if and only if
$r^{2}\mid n$.  Consequently, writing $n=m^{2}s$ with $s$ squarefree, the maximal
Tannakian subcategory of $\mathcal{Z}(\mathcal{N}_G)_{\mathrm{pt}}$ is $\mathcal{C}(H_m,q|_{H_m})\simeq\mathrm{Rep}(C_m)$ and this depends only on $n$, not on the quadratic form $q$ nor on $k$.
\end{proposition}
 
\begin{proof}
The elements of $H_r$ are $jn/r$ for $j\in\mathbb{Z}$, and $q(jn/r)=\kappa^{j^2n^2/r^2}$. As $\kappa$ has order $n$, $q|_{H_r}\equiv1$ if and only if $n\mid j^{2}n^{2}/r^{2}$, that is, $r^{2}\mid j^{2}n$, for all $j$; taking $j=1$, this holds if and only if $r^{2}\mid n$, which gives the first claim.
 
For a prime $p$ we have $2v_p(r)\le v_p(n)$ if and only if
$v_p(r)\le\lfloor v_p(n)/2\rfloor=v_p(m)$, so the divisors $r$ of $n$ with
$r^{2}\mid n$ are exactly the divisors of $m$, and $H_r\subseteq H_m$ for every $r\mid m$. So $H_m$ is the unique maximal one, and $\mathcal{C}(H_m,q|_{H_m})\simeq\mathrm{Rep}(\widehat{H_m})\simeq\mathrm{Rep}(C_m)$.  Finally the criterion $r^{2}\mid n$ does
not involve $\kappa$ or $k$.
\end{proof}

\par The cocycles in Corollary \ref{cor:lag} and Proposition \ref{prop:Dpt} below are detected by self-braidings. Identify $C_m$ with $\{0,1,\dots,m-1\}$, and write $[v]$ for the residue of $v$ modulo $m$. The classes in $H^{3}(C_m,\mathbb{C}^{\times})\cong\mathbb{Z}/m\mathbb{Z}$ are those of the cocycles
\[
\omega_a(u,v,v')=\exp\bigl(2\pi i\,a\,u\,(v+v'-[v+v'])/m^2\bigr),\qquad a\in\{0,1,\dots,m-1\},
\]
and the class of $\omega_a$ has order $m/\gcd(a,m)$. A half-braiding on the simple object $u$ of $\mathrm{Vec}^{\omega_a}_{C_m}$ is given by scalars $\chi(v)$ with $\chi(v)\chi(v')=\chi([v+v'])\exp\bigl(\pm2\pi i\,a\,u\,(v+v'-[v+v'])/m^2\bigr)$, where the sign depends only on conventions. The solutions are $\chi(v)=\exp(\pm2\pi i\,a\,u\,v/m^2+2\pi i\,\ell v/m)$ with $\ell\in\mathbb{Z}/m\mathbb{Z}$, and the self-braiding of the resulting invertible object $Y$ is $\chi(u)$. Hence
\begin{equation}\label{eq:twistm}
c_{Y,Y}^{\,m}=\exp(\pm2\pi i\,a\,u^2/m)
\end{equation}
for every simple object $Y$ of $\mathcal{Z}(\mathrm{Vec}^{\omega_a}_{C_m})$ lying over $u$. In particular, if $c_{Y,Y}^{\,m}$ is a primitive $m$-th root of unity for some $Y$, then $\gcd(a,m)=1$, so $\omega_a$ has order $m$.

\begin{corollary}\label{cor:lag}
Assume $n>1$ is odd.  The maximal Tannakian subcategory of
$\mathcal{Z}(\mathcal{N}_G)_{\mathrm{pt}}$ is Lagrangian if and only if $n$ is a
perfect square, say $n=m^{2}$.  In that case
$\mathcal{Z}(\mathcal{N}_G)_{\mathrm{pt}}\simeq
\mathcal{Z}(\mathrm{Vec}^{\omega}_{C_m})$ for some
$\omega\in H^{3}(C_m,\mathbb{C}^{\times})\cong\mathbb{Z}/m\mathbb{Z}$ of order $m$, and
$[\mathcal{Z}(\mathcal{N}_G)_{\mathrm{pt}}]
=[\mathcal{Z}(\mathcal{N}_G)_{\mathrm{ad}}]=0$ in the Witt group of
nondegenerate braided fusion categories.
\end{corollary}
 
\begin{proof}
For $n$ odd, $\mathcal{Z}(\mathcal{N}_G)_{\mathrm{pt}}\simeq\mathcal{C}(G,q)$ is nondegenerate by
Proposition \ref{cor:cyclic}, hence modular, and its dimension is $n$; so
Lagrangian means $m^{2}=n$, i.e.\ $s=1$.  A modular category containing a
Lagrangian subcategory $\mathrm{Rep}(C_m)$ is braided equivalent to
$\mathcal{Z}(\mathrm{Vec}^{\omega}_{C_m})$ for some $\omega$ \cite[Theorem 4.64]{MR2609644}. To see that $\omega$ has order $m$, note that the object of $\mathcal{Z}(\mathrm{Vec}^{\omega}_{C_m})$ corresponding to $Y_1$ has self-braiding $q(1)=\kappa$. As $\kappa$ has order $n=m^2$, $\kappa^m$ is a primitive $m$-th root of unity, so $\omega$ has order $m$ by \eqref{eq:twistm}.  A Drinfeld center is Witt trivial \cite{MR3039775}, and
$[\mathcal{Z}(\mathcal{N}_G)_{\mathrm{ad}}]
=-[\mathcal{Z}(\mathcal{N}_G)_{\mathrm{pt}}]$ because
$\mathcal{Z}(\mathcal{N}_G)_{\mathrm{pt}}\boxtimes
\mathcal{Z}(\mathcal{N}_G)_{\mathrm{ad}}\simeq\mathcal{Z}(\mathcal{N}_G)$ is Witt
trivial; this factorization holds by Corollary \ref{cor:oddfactor}.
\end{proof}


\subsection{The associated de-equivariantization}\label{subsec:assocdeeq}
 
In this subsection we write $G$ and its quotients multiplicatively, with identity $e$. Let $r^{2}\mid n$ and put
$\mathcal{E}:=\mathcal{C}(H_r,q|_{H_r})\simeq\mathrm{Rep}(C_r)$.  By
\eqref{eq:star} the composite
$\mathcal{E}\hookrightarrow\mathcal{Z}(\mathcal{N}_G)\xrightarrow{F}\mathcal{N}_G$
sends the $r$ simple objects $Y_h$, $h\in H_r$, to the $r$ distinct invertible
objects $h\in G$, hence is fully faithful.  Therefore $\mathcal{N}_G$ is a
$C_r$-equivariantization.  Precisely, let $A_r$ denote the regular algebra of
$\mathcal{E}$, a connected \'etale algebra in $\mathcal{Z}(\mathcal{N}_G)$, and put
$A:=F(A_r)=\bigoplus_{h\in H_r}h$, an algebra in $\mathcal{N}_G$. We define
\[
\mathcal{D}:=(\mathcal{N}_G)_{C_r}:=\mathrm{Rep}_{\mathcal{N}_G}(A_r)
\]
as in \cite[Definition 3.17]{MR3039775}: the category of right $A$-modules in $\mathcal{N}_G$, with the tensor product $\otimes_A$ for which a right $A$-module $M$ is an $A$-bimodule via $A\otimes M\xrightarrow{\gamma^{-1}}M\otimes A\to M$, where $\gamma\colon M\otimes A\to A\otimes M$ is the half-braiding of $A_r$. This is the de-equivariantization of $\mathcal{N}_G$ by $\mathcal{E}$ in the sense of \cite[Section 4.2.4]{MR2609644}. There left $A$-modules are used, made into bimodules through $\gamma$; composing a right action with $\gamma^{-1}$ identifies right with left $A$-modules, and the two resulting bimodule structures coincide, so the tensor product \cite[(68)]{MR2609644} agrees with $\otimes_A$. In particular:
\begin{itemize}
\item $\mathcal{D}$ is a fusion category with $\operatorname{FPdim}(\mathcal{D})=\operatorname{FPdim}(\mathcal{N}_G)/r$ \cite[Remark 4.21 and Proposition 4.26]{MR2609644}, and the free module functor $X\mapsto X\otimes A$ is a surjective tensor functor $\mathcal{N}_G\to\mathcal{D}$ \cite[Section 4.2.6]{MR2609644}, \cite[Remark 3.19]{MR3039775};
\item the group $\Pi\cong C_r$ of automorphisms of the algebra $A_r$ acts on $\mathcal{D}$ by tensor autoequivalences, and $\mathcal{N}_G\simeq\mathcal{D}^{\Pi}$, with the free module functor corresponding to the forgetful functor $\mathcal{D}^\Pi\to\mathcal{D}$ \cite[Proposition 4.19 and Lemma 4.24]{MR2609644};
\item there is a braided tensor equivalence \cite[Theorem 3.20]{MR3039775}
\begin{equation}\label{eq:ZD}
  \mathcal{Z}(\mathcal{D})\simeq\mathcal{Z}(\mathcal{N}_G)^0_{A_r}.
\end{equation}
\end{itemize}
Nothing below uses the parity of $n$.

\begin{proposition}\label{prop:deeq}
With the notation above, $\mathcal{D}$ has exactly $n/r$ invertible objects,
forming the group $G/H_r\cong C_{n/r}$, and exactly $r$ noninvertible simple objects $M_1,\dots,M_r$. These satisfy $x\otimes A\cong M_1\oplus\cdots\oplus M_r$, they are permuted freely and transitively by $\Pi$, and each has Frobenius--Perron dimension $\delta:=d/r$, where
\begin{equation}\label{eq:delta}
  \delta^{2}=\frac{kn}{r}\,\delta+\frac{n}{r^{2}} .
\end{equation}
\end{proposition}

\begin{proof}
For $g\in G$ we have $g\otimes A=\bigoplus_{h\in H_r}gh$, the coset $gH_r$, and
$\operatorname{End}_A(g\otimes A)\cong
\operatorname{Hom}_{\mathcal{N}_G}(g,g\otimes A)\cong\mathbb{C}$, so $g\otimes A$
is a simple $A$-module of dimension $1$ in $\mathcal{D}$; two such agree
precisely when the cosets agree, giving the group $G/H_r\cong C_{n/r}$. Moreover $\mathrm{Hom}_A(x\otimes A,g\otimes A)\cong\mathrm{Hom}_{\mathcal{N}_G}(x,g\otimes A)=0$. As the free module functor is surjective, the simple objects of $\mathcal{D}$ are therefore the $g\otimes A$ and the simple summands of $x\otimes A$, and no $g\otimes A$ is a summand of $x\otimes A$.

Under $\mathcal{N}_G\simeq\mathcal{D}^\Pi$, the objects $g\otimes A$ and $x\otimes A$ are forgetful images of equivariant objects, so each is isomorphic to its image under every element of $\Pi$. Hence $\Pi$ fixes each simple object $g\otimes A$ and permutes the simple summands of $x\otimes A$. Recall that the simple objects of $\mathcal{D}^\Pi$ lying over a $\Pi$-orbit $\Omega$ of simple objects of $\mathcal{D}$ correspond to the irreducible projective representations of the stabilizer $\Pi_M\subseteq\Pi$ of an object $M\in\Omega$, for a $2$-cocycle on $\Pi_M$ determined by the action, and that the forgetful functor sends each of these simple objects to $\bigoplus_{M'\in\Omega}M'$ with multiplicity the dimension of the corresponding representation \cite[Corollary 2.13]{MR3059899}. As $\Pi_M$ is a subgroup of the cyclic group $\Pi$, $H^2(\Pi_M,\mathbb{C}^\times)=0$, so these representations are the $|\Pi_M|$ linear characters of $\Pi_M$, each of dimension $1$. Now let $\Omega$ be the $\Pi$-orbit of a simple summand $M$ of $x\otimes A$. Then $\Omega$ consists of summands of $x\otimes A$, so it contains no $g\otimes A$, and $x$ is the only simple object of $\mathcal{N}_G\simeq\mathcal{D}^\Pi$ lying over $\Omega$. Hence $|\Pi_M|=1$, and $x\otimes A\cong\bigoplus_{M'\in\Omega}M'$ with $|\Omega|=|\Pi|=r$. That is, $x\otimes A$ is the direct sum of $r$ pairwise non-isomorphic simple objects $M_1,\dots,M_r$, permuted freely and transitively by $\Pi$, and these are all the simple objects of $\mathcal{D}$ besides the $g\otimes A$.

Tensor autoequivalences preserve Frobenius--Perron dimensions, so the $M_i$ all have the same dimension, and $r\operatorname{FPdim}(M_i)=\operatorname{FPdim}(x)=d$ because the free module functor is a tensor functor. Thus $\operatorname{FPdim}(M_i)=\delta=d/r\notin\mathbb{Q}$; in particular no $M_i$ is invertible. Dividing $d^{2}=knd+n$ \eqref{eq:ids} by $r^{2}$ gives \eqref{eq:delta}.
\end{proof}
 
\begin{definition}\label{def:haagerup}
Let $n,k\in\mathbb{Z}_{\geq1}$ and let $z\in C_n$ with $z^2=e$. The \emph{stretched Haagerup--Izumi fusion ring} $H_k^z(C_n)$ is the fusion ring with basis $\{g,\ g\rho:g\in C_n\}$, and nontrivial fusion rules
\begin{equation}\label{eq:haagring}
  \rho g=g^{-1}\rho,\qquad\text{and}\qquad \rho^2=z+k\sum_{g\in C_n}g\rho .
\end{equation}
Its duality is given by $g^\ast=g^{-1}$ and $(g\rho)^\ast=zg\rho$, and $(g\rho)(h\rho)=zgh^{-1}+k\sum_{f\in C_n}f\rho$ for $g,h\in C_n$. Each $g\rho$ has Frobenius--Perron dimension equal to the positive root of $t^2=knt+1$. We write $H_k(C_n):=H_k^e(C_n)$, in which every $g\rho$ is self-dual; if $z\neq e$, then $n$ is even, $z$ is the element of order $2$, and no $g\rho$ is self-dual. For $k=1$ we call $H_1^z(C_n)$ a \emph{Haagerup--Izumi fusion ring}; the rings $H_1(C_n)$ are the usual Haagerup--Izumi fusion rings (see e.g.\ \cite{MR2837122}), and $H_1(C_3)$ is the fusion ring of the principal even part of the Haagerup subfactor.
\end{definition}

\begin{theorem}\label{thm:haagerup}
Let $m>1$, and let $\mathcal{N}_G$ be a near-group fusion category of type $C_{m^2}+km^2$ satisfying \eqref{eq:star}, for instance $k=1$. Then the maximal Tannakian subcategory of $\mathcal{Z}(\mathcal{N}_G)_\mathrm{pt}$ is $\mathrm{Rep}(C_m)$, and the de-equivariantization $\mathcal{D}=(\mathcal{N}_G)_{C_m}$ satisfies $K(\mathcal{D})\cong H_k(C_m)$ if $m$ is odd and $K(\mathcal{D})\cong H_k^z(C_m)$ if $m$ is even, where $z\in C_m$ is the element of order $2$.
\end{theorem}

\begin{proof}
Here $n=m^{2}$, so Proposition \ref{prop:maxtann} gives the first
claim.  Apply Proposition \ref{prop:deeq} with $r=m$:
$\mathcal{D}$ has $n/m=m$ invertibles, forming $C_m$, and $m$ noninvertible simples of Frobenius--Perron dimension
$\delta=d/m$ with $\delta^{2}=(kn/m)\delta+n/m^{2}=km\delta+1$ by
\eqref{eq:delta}.  Note that $\delta\notin\mathbb{Q}$, as $d\notin\mathbb{Q}$.

The group $C_m=\mathcal{O}(\mathcal{D}_\mathrm{pt})$ acts on the set of $m$ noninvertible
simples by left tensor multiplication. Let $M$ be a noninvertible simple object with stabilizer $\mathrm{Stab}(M)\subseteq C_m$.  Then $[g,M\otimes M^{*}]=[M,g\otimes M]$ is $1$ for $g\in\mathrm{Stab}(M)$ and $0$ otherwise, so
$\delta^{2}=|\mathrm{Stab}(M)|+\sum_ia_i\delta$ for some $a_i\in\mathbb{Z}_{\ge0}$.  Comparing
with $\delta^{2}=km\delta+1$ gives $(km-\sum_ia_i)\delta=|\mathrm{Stab}(M)|-1$, so $|\mathrm{Stab}(M)|=1$ by
irrationality of $\delta$.  Hence the action is free, and therefore transitive, as there are $m$ noninvertible simple objects.  The same argument applied to $M^\ast\otimes M$ shows that
right tensor multiplication by $C_m$ is also free and transitive.  Fix a
noninvertible simple object $M$.  Then every noninvertible simple object is
isomorphic to $h\otimes M$ for a unique $h\in C_m$, and there is a bijection
$\vartheta:C_m\to C_m$ with $M\otimes g\cong\vartheta(g)\otimes M$ for all $g\in C_m$.

Consequently $[g,M\otimes M^{*}]=0$ for $g\neq e$ while
$[\mathbbm{1},M\otimes M^{*}]=1$, so
\[
  M\otimes M^{*}\cong\mathbbm{1}\oplus\bigoplus_{h\in C_m}a_h\,(h\otimes M)
\]
with $\sum_ha_h\delta=\delta^{2}-1=km\delta$, i.e.\ $\sum_ha_h=km$.  To see that
the $a_h$ are all equal, let $\Pi\cong C_m$ be the group acting on $\mathcal{D}$
by tensor autoequivalences with $\mathcal{N}_G\simeq\mathcal{D}^{\Pi}$, as in Section \ref{subsec:assocdeeq}.  By Proposition \ref{prop:deeq} and its proof, $\Pi$ fixes each invertible object of $\mathcal{D}$ and permutes the noninvertible
simple objects freely and transitively, so for each $\phi\in\Pi$ we have
$\phi(M)\cong g_\phi\otimes M$ for a unique $g_\phi\in C_m$, and
$\phi\mapsto g_\phi$ is a bijection $\Pi\to C_m$.  Applying $\phi$ to the
decomposition of $M\otimes M^\ast$, and using $\phi(h)\cong h$, gives
\[
  \phi(M\otimes M^{*})\cong\mathbbm{1}\oplus\bigoplus_{h\in C_m}a_h\,(hg_\phi\otimes M),
\]
while
\[
  \phi(M)\otimes\phi(M)^{*}\cong g_\phi\otimes(M\otimes M^{*})\otimes g_\phi^{-1}
  \cong\mathbbm{1}\oplus\bigoplus_{h\in C_m}a_h\,\bigl(g_\phi h\,\vartheta(g_\phi^{-1})\otimes M\bigr).
\]
Comparing multiplicities gives $a_u=a_{u\vartheta(g_\phi^{-1})^{-1}}$ for all
$u\in C_m$.  As $\phi$ runs over $\Pi$, $g_\phi^{-1}$ runs over $C_m$, hence so
does $\vartheta(g_\phi^{-1})^{-1}$; therefore $a_h$ is independent of $h$, and
$a_h=k$ for all $h\in C_m$. That is,
\begin{equation}\label{eq:haag}
  M\otimes M^{*}\cong\mathbbm{1}\oplus k\!\!\bigoplus_{g\in C_m}\!\!(g\otimes M).
\end{equation}

We now show that $\vartheta(g)=g^{-1}$ for all $g\in C_m$. Write $M^\ast\cong u\otimes M$ with $u\in C_m$. Tensoring $M\otimes g\cong\vartheta(g)\otimes M$ on the right with $g^{-1}$ gives $M\otimes g^{-1}\cong\vartheta(g)^{-1}\otimes M$, that is, $\vartheta(g^{-1})=\vartheta(g)^{-1}$. Let $\phi\in\Pi$. As $\phi$ is a tensor autoequivalence fixing the invertible object $u$,
\[
  ug_\phi\otimes M\cong\phi(u\otimes M)\cong\phi(M^\ast)\cong\phi(M)^\ast\cong M^\ast\otimes g_\phi^{-1}\cong u\otimes M\otimes g_\phi^{-1}\cong u\,\vartheta(g_\phi)^{-1}\otimes M .
\]
Since the left action is free, $\vartheta(g_\phi)=g_\phi^{-1}$, and since $g_\phi$ runs over $C_m$, $\vartheta(g)=g^{-1}$ for all $g\in C_m$.

Finally, $M\cong M^{\ast\ast}\cong(u\otimes M)^\ast\cong M^\ast\otimes u^{-1}\cong u\otimes M\otimes u^{-1}\cong u^2\otimes M$, so $u^2=e$. For every $h\in C_m$ we have $(h\otimes M)^\ast\cong M^\ast\otimes h^{-1}\cong u\otimes M\otimes h^{-1}\cong u\otimes(h\otimes M)$, so $X^\ast\cong u\otimes X$ for every noninvertible simple $X$. As $M\cong u\otimes M^\ast$ and $M\otimes u\cong u\otimes M$, \eqref{eq:haag} gives
\[
  M\otimes M\cong u\otimes(M\otimes M^\ast)\cong u\oplus k\!\!\bigoplus_{g\in C_m}\!\!(g\otimes M),
\]
while $(h\otimes M)\otimes g\cong hg^{-1}\otimes M$ and $(h\otimes M)\otimes(h'\otimes M)\cong hh'^{-1}\otimes(M\otimes M)$. These are the products and the duality of $H_k^u(C_m)$ under $g\mapsto g$, $g\otimes M\mapsto g\rho$. If $m$ is odd, then $u^2=e$ forces $u=e$.

Now let $m$ be even and let $z\in C_m$ be the element of order $2$. As $u^2=e$, we have $u\in\{e,z\}$, and we show that $u=z$. The category $\mathcal{D}$ is pseudounitary, since $\dim(\mathcal{D})=\dim(\mathcal{N}_G)/m=\mathrm{FPdim}(\mathcal{N}_G)/m=\mathrm{FPdim}(\mathcal{D})$ \cite[Proposition 4.26]{MR2609644}, and we equip it with its canonical spherical structure \cite[Proposition 9.5.1]{MR3242743}. Let $w\in G$ correspond to $m/2\in\mathbb{Z}/n\mathbb{Z}$. Then $w\notin H_m$ and $w^2\in H_m$, so the invertible object $w\otimes A$ of $\mathcal{D}$ is $z$ (Proposition \ref{prop:deeq}). As $4\mid n$, $q(w)=\kappa^{m^2/4}=\kappa^{n/4}$ is a primitive fourth root of unity, and $b(w,h)=\kappa^{2\cdot(m/2)\cdot jm}=\kappa^{jn}=1$ for each $h\in H_m$, say $h$ corresponding to $jm$. The first identity says that $c_{A_m,Y_w}\circ c_{Y_w,A_m}=\mathrm{id}$, as $A_m=\bigoplus_{h\in H_m}Y_h$. Let $\widehat{Z}:=Y_w\otimes A_m$ be the free $A_m$-module, with action $\mathrm{id}_{Y_w}\otimes\mu$, where $\mu$ is the multiplication of $A_m$. By the hexagon axioms, $c_{A_m,\widehat{Z}}\circ c_{\widehat{Z},A_m}=(\mathrm{id}_{Y_w}\otimes c_{A_m,A_m})\circ(c_{A_m,Y_w}c_{Y_w,A_m}\otimes\mathrm{id}_{A_m})\circ(\mathrm{id}_{Y_w}\otimes c_{A_m,A_m})=\mathrm{id}_{Y_w}\otimes c_{A_m,A_m}^2$, and $\mu\circ c_{A_m,A_m}=\mu$ as $A_m$ is commutative; so $\widehat{Z}$ is a local $A_m$-module (Section \ref{app:correspond}). It has dimension $\dim(Y_w\otimes A_m)/\dim(A_m)=1$, and $\widehat{Z}\otimes_{A_m}\widehat{Z}\cong Y_{w^2}\otimes A_m\cong A_m$. The simple summands of its underlying object are the $Y_{wh}$ with $h\in H_m$, each of twist $q(wh)=q(w)q(h)b(w,h)=q(w)$. As $q|_{H_m}\equiv1$, we have $\theta_{A_m}=\mathrm{id}$, so by Proposition \ref{twist}(1) the twist of $\widehat{Z}$ in the ribbon structure inherited from $\mathcal{Z}(\mathcal{N}_G)$ is $q(w)$. The equivalence \eqref{eq:ZD} lifts the forgetful functor $\mathcal{Z}(\mathcal{N}_G)^0_{A_m}\to\mathcal{D}$ \cite[(22) and Theorem 3.20]{MR3039775}, so $\widehat{Z}$ is a self-dual invertible object of $\mathcal{Z}(\mathcal{D})$ lying over $z$. Both the inherited ribbon structure and the canonical ribbon structure of $\mathcal{Z}(\mathcal{D})$ have positive dimensions, so they agree under \eqref{eq:ZD} (compare the argument after \eqref{eq:S0}). Hence $\widehat{Z}$ has twist $q(w)$ in $\mathcal{Z}(\mathcal{D})$, and as $\dim(\widehat{Z})=1$, its self-braiding is also $q(w)$. Its double braiding with itself is $q(w)^2=-1$, so the fusion subcategory of $\mathcal{Z}(\mathcal{D})$ generated by $\widehat{Z}$ is a modular category of rank $2$ with dimensions $1,1$ and twists $1,q(w)$. The formula for Frobenius--Schur indicators in modular categories \cite[Theorem 7.5]{MR2313527} gives
\[
\nu_2(\widehat{Z})=\tfrac12\bigl(q(w)^{-2}+q(w)^{2}\bigr)=-1,
\]
and hence $\nu_2(z)=-1$, as the forgetful functor $\mathcal{Z}(\mathcal{D})\to\mathcal{D}$ is a pivotal tensor functor (Section \ref{subsec:nucentral}).

Suppose that $u=e$. Then every noninvertible simple object of $\mathcal{D}$ is self-dual, and Lemma \ref{lem:nucentral}, applied to $\widehat{Z}$ and $M$, gives $\nu_2(z\otimes M)=\nu_2(z)\nu_2(M)=-\nu_2(M)$. On the other hand, let $\phi\in\Pi$ with $g_\phi=z$, so that $\phi(M)\cong z\otimes M$. The tensor autoequivalence $\phi$ preserves Frobenius--Perron dimensions, hence preserves the canonical spherical structure of $\mathcal{D}$ \cite[Proposition 9.5.1]{MR3242743}, and therefore preserves Frobenius--Schur indicators (Section \ref{subsec:nucentral}; see also \cite[Section 6]{MR2381536}). Thus $\nu_2(z\otimes M)=\nu_2(M)$, so $\nu_2(M)=0$, which is impossible for a self-dual simple object. Hence $u=z$.
\end{proof}

\begin{proposition}\label{prop:Dpt}
Let $n=m^2$ with $m>1$, and let $\mathcal{D}=(\mathcal{N}_G)_{C_m}$ be as in Theorem \ref{thm:haagerup}. Then $\mathcal{D}_\mathrm{pt}\simeq\mathrm{Vec}^{\omega'}_{C_m}$ for some $\omega'\in H^3(C_m,\mathbb{C}^\times)$ of order $m$.
\end{proposition}

\begin{proof}
Restricting half-braidings to $(\mathcal{N}_G)_\mathrm{pt}=\mathrm{Vec}_G$ gives a braided tensor functor $R\colon\mathcal{Z}(\mathcal{N}_G)_\mathrm{pt}\to\mathcal{Z}(\mathrm{Vec}_G)$. In the proof of Proposition \ref{correspond}, this restriction corresponds to the free module functor to $\mathcal{Z}(\mathcal{N}_G)_{A_G}$, so $R(Y_g)=(g,t_g)$ for every $g\in G$ by Lemma \ref{lem:liftable}(2). Hence $A_0:=R(A_m)=\bigoplus_{h\in H_m}(h,t_h)$ is a connected \'etale algebra in $\mathcal{Z}(\mathrm{Vec}_G)$ with underlying object $A$. The simple objects of $\mathcal{D}_\mathrm{pt}$ are the $A$-modules $g\otimes A$ (Proposition \ref{prop:deeq}), and their tensor products involve only the half-braidings of $A_m$ with objects of $\mathrm{Vec}_G$. So $\mathcal{D}_\mathrm{pt}=\mathrm{Rep}_{\mathrm{Vec}_G}(A_0)$, and $\mathcal{Z}(\mathcal{D}_\mathrm{pt})\simeq\mathcal{Z}(\mathrm{Vec}_G)^0_{A_0}$ \cite[Theorem 3.20]{MR3039775}.

Here $t_g(h)=\kappa^{gh}$ for all $g,h\in\mathbb{Z}/n\mathbb{Z}$; in particular $t_g(h)=t_h(g)$. Let $w\in G$ correspond to $1$, and put $W:=(w,t_w^{-1})\in\mathcal{Z}(\mathrm{Vec}_G)$. By \eqref{eq:S0}, the double braiding of $W$ with $(h,t_h)$ is $t_w(h)^{-1}t_h(w)=1$. So, by the hexagon argument in the proof of Theorem \ref{thm:haagerup}, the free module $V:=W\otimes A_0$ is a local $A_0$-module of dimension $1$, that is, an invertible object of $\mathcal{Z}(\mathcal{D}_\mathrm{pt})$. The simple summands of $V$ are the invertible objects $W\otimes(h,t_h)$, $h\in H_m$, each of twist $(t_w^{-1}t_h)(wh)=t_w(w)^{-1}\,t_h(h)=\kappa^{-1}$, using $t_w(h)=t_h(w)$ and $q|_{H_m}\equiv1$. The latter also gives $\theta_{A_0}=\mathrm{id}$, so by Proposition \ref{twist}(1) the twist of $V$ in the ribbon structure inherited from $\mathcal{Z}(\mathrm{Vec}_G)$ is $\kappa^{-1}$, and as $\dim(V)=1$ its self-braiding is also $\kappa^{-1}$. As $\mathcal{D}_\mathrm{pt}$ is pointed with group $G/H_m\cong C_m$, we have $\mathcal{D}_\mathrm{pt}\simeq\mathrm{Vec}^{\omega'}_{C_m}$ for some $\omega'$. Since $\kappa$ has order $n=m^2$, $\kappa^{-m}$ is a primitive $m$-th root of unity, so $\omega'$ has order $m$ by \eqref{eq:twistm}.
\end{proof}

\begin{note}\label{rem:evenm}
For $m=2$ and $k=1$, the ring $H_1^z(C_2)$ is categorified by the even part of a subfactor \cite[Remark 1.2]{MR3229513}. In the unitary setting, Izumi de-equivariantizes the two $C^\ast$ near-group categories for $C_4$ by $C_2$, which is not Lagrangian in his sense, and obtains the $1$-supertransitive subfactors of Liu, Morrison and Penneys \cite[Example 12.12]{MR3635673}, in accordance with Theorem \ref{thm:haagerup}. By contrast, for $G=C_2\times C_2$ the de-equivariantization by a subgroup of order $2$ is the even part of $SU(2)_6$ \cite[Example 12.11]{MR3635673}, whose fusion ring is $H_1(C_2)$. For even $m$ the proof of Theorem \ref{thm:haagerup} uses that $G$ is cyclic, through the fact that $q(w)$ is a primitive fourth root of unity. It also shows that $\nu_2(z)=-1$. For $k=2$, a near-group fusion category of type $C_4+8$ satisfying \eqref{eq:star} would categorify $H_2^z(C_2)$, the only ring in \cite[Theorem 1.1]{MR3229513} without a known categorification; there is no $C^\ast$ near-group category of this type \cite[Theorem 10.22]{MR3635673}.
\end{note}

\begin{note}\label{rem:morita}
Let $n=m^2$ with $m>1$ odd. Then $\mathcal{N}_G$ is $G$-factorizable (Corollary \ref{cor:oddfactor}), and $\mathcal{Z}(\mathcal{N}_G)_\mathrm{pt}\simeq\mathcal{Z}(\mathrm{Vec}^\omega_{C_m})$ with $\mathcal{E}\simeq\mathrm{Rep}(C_m)$ Lagrangian, where $\omega$ is as in Corollary \ref{cor:lag}. As $A_m\in\mathcal{Z}(\mathcal{N}_G)_\mathrm{pt}$ and $(\mathcal{Z}(\mathcal{N}_G)_\mathrm{pt})^0_{A_m}\simeq\mathrm{Vec}$, the free module functor $W\mapsto W\otimes A_m$ is a braided equivalence $\mathcal{Z}(\mathcal{N}_G)_\mathrm{ad}\simeq\mathcal{Z}(\mathcal{N}_G)^0_{A_m}$ (as in \cite[Lemma 3.29]{MR3039775}; see \cite[Remark 4.7]{MR3039775}), so $\mathcal{Z}(\mathcal{D})\simeq\mathcal{Z}(\mathcal{N}_G)_\mathrm{ad}$ by \eqref{eq:ZD}. Thus Section \ref{sec:factorizable} describes $\mathcal{Z}(\mathcal{D})$, the double of the stretched Haagerup--Izumi category $\mathcal{D}$, as a modular category: its simple objects, dimensions, and $S$- and $T$-matrices. Its forgetful functor to $\mathcal{D}$ is the composite $\mathcal{Z}(\mathcal{N}_G)_\mathrm{ad}\to\mathcal{Z}(\mathcal{N}_G)^0_{A_m}\to\mathcal{D}$. Moreover $\mathcal{Z}(\mathrm{Vec}^\omega_{C_m}\boxtimes\mathcal{D})\simeq\mathcal{Z}(\mathcal{N}_G)_\mathrm{pt}\boxtimes\mathcal{Z}(\mathcal{N}_G)_\mathrm{ad}\simeq\mathcal{Z}(\mathcal{N}_G)$, so $\mathcal{N}_G$ is Morita equivalent to $\mathrm{Vec}^\omega_{C_m}\boxtimes\mathcal{D}$ \cite[Theorem 3.1]{MR2735754}. This is a converse, valid for every $k$ under \eqref{eq:star}, to Izumi's construction of near-group categories from categories of Haagerup type in the unitary setting with $k=1$ \cite[Theorem 12.9 and Remark 12.10]{MR3635673}.

The same arguments apply when $G$ is any abelian group of odd order $m^2$ satisfying \eqref{eq:star} and $(G,q)$ contains a subgroup $L$ of order $m$ with $q|_L\equiv1$: then $\mathcal{N}_G$ is Morita equivalent to $\mathrm{Vec}^{\omega}_{L}\boxtimes\mathcal{D}$ for some $\omega\in H^3(L,\mathbb{C}^\times)$, where $\mathcal{D}$ is the de-equivariantization by $L$, defined as in Section \ref{subsec:assocdeeq}, and if $L$ and $G/L$ are cyclic, for instance if $m=p$ is prime, then $\mathcal{D}$ categorifies $H_k(C_m)$. This applies to $G=C_p\times C_p$ whenever $q$ is isotropic, and then $\omega$ is trivial: the values of $q$ lie in $\mu_p$ (Proposition \ref{prop:block}(1)), so $c_{Y,Y}^{\,p}=1$ for every simple object $Y$ of $\mathcal{Z}(\mathrm{Vec}^{\omega_a}_{C_p})$, where $\omega$ is the class of $\omega_a$, and \eqref{eq:twistm} with $u\neq0$ gives $a=0$. For instance, Izumi de-equivariantizes the $C^\ast$ near-group category for $C_3\times C_3$ by two subgroups $L$ of order $3$ and obtains two categorifications of the Haagerup fusion ring \cite[Example 12.13]{MR3635673}; such a de-equivariantization requires $\mathcal{C}(L,q|_L)$ to be Tannakian (Corollary \ref{cor:k1}), so $q|_L\equiv1$ \cite[Example 2.48]{MR2609644}.
\end{note}


\section{Applications to arbitrary categorifications}\label{sec:application}

In this section $\mathcal{N}_G$ is a near-group fusion category of type $G+kn$, where $n=|G|$ and $k\in\mathbb{Z}_{\geq1}$, with the conventions of Section \ref{sec:neargroup}; neither \eqref{eq:star} nor $G$-factorizability is assumed. We use $d$, $\Delta=n(k^2n+4)$, $d'$, $N$ and $\mathfrak{d}$ as in Section \ref{subsec:indicators}, and write $\Delta=c^2\Delta_0$ with $c\in\mathbb{Z}_{\geq1}$ and $\Delta_0$ squarefree. As $\Delta$ is not a square (beginning of Section \ref{sec:neargroup}), $\Delta_0>1$, and the conductor $\mathfrak{f}=|\mathfrak{d}|$ of $\mathbb{Q}(\sqrt\Delta)=\mathbb{Q}(\sqrt{\Delta_0})$ is $\Delta_0$ if $\Delta_0\equiv1\pmod4$ and $4\Delta_0$ otherwise.

\par Write $R(G,\ell)$ for the near-group fusion ring with basis $G\cup\{x\}$ and $x^2=\sum_{g\in G}g+\ell x$, so that $R(G,kn)$ is $G+kn$. Apart from one well-understood family, these are the only near-group rings with $\ell>0$ that can be categorified.

\begin{lemma}\label{lem:levels}
Let $G$ be a finite abelian group of order $n$, and let $\ell>0$. If $R(G,\ell)$ is categorifiable, then either $n\mid\ell$, or $\ell=n-1$, $G$ is cyclic and $n+1$ is a prime power.
\end{lemma}

\begin{proof}
If $0<\ell<n-1$, then a categorifiable $R(G,\ell)$ would be noncommutative \cite[Proposition 3.2]{MR4655273}, whereas $R(G,\ell)$ is commutative as $G$ is abelian. So if $\ell<n$, then $\ell=n-1$, and then $G$ is cyclic with $n+1$ a prime power \cite[Theorem 1.2]{Siehler03}. If $\ell\geq n$, then $\ell^2<\ell^2+4n<(\ell+2)^2$, and $\ell^2+4n\neq(\ell+1)^2$ by parity. So $\frac12(\ell+\sqrt{\ell^2+4n})\notin\mathbb{Q}$, and $n\mid\ell$ by \cite[Proposition 3.3]{MR4655273}.
\end{proof}

\subsection{A conductor bound}

\begin{proposition}\label{prop:conductor}
Let $G$ be a finite abelian group of order $n$ and $k\in\mathbb{Z}_{\geq1}$. If the stretched near-group fusion ring $G+kn$ is categorifiable, then
\begin{equation}\label{eq:condbound}
    k^2\mathfrak{f}\leq\Delta=n(k^2n+4).
\end{equation}
Equivalently, $c\geq k$ if $\Delta_0\equiv1\pmod4$, and $c\geq2k$ otherwise.
\end{proposition}

\begin{proof}
In the notation of \cite{schopieray2026dimensionblocksphericalfusion}, $\mathbb{K}_1=\mathbb{Q}(d)=\mathbb{Q}(\sqrt\Delta)$ and $G_1=\{1,\sigma\}$, where $\sigma$ is the nontrivial automorphism of $\mathbb{Q}(\sqrt\Delta)$. We have $\dim(\mathcal{N}_G)=n(2+kd)$ and $\sigma(\dim(\mathcal{N}_G))=n(2+kd')$. By \eqref{eq:ddprime}, their difference is $nk\sqrt\Delta$ and their product is $n^2(4+k^2n)=n\Delta$, so the nontrivial character $\chi$ of $G_1$ has $|\hat g_1(\chi)|=\bigl|\dim(\mathcal{N}_G)^{-1}-\sigma(\dim(\mathcal{N}_G))^{-1}\bigr|=k/\sqrt\Delta$. Then \cite[Theorem 5.1]{schopieray2026dimensionblocksphericalfusion} gives \eqref{eq:condbound}. The second form follows from $\Delta=c^2\Delta_0$ and the description of $\mathfrak{f}$ above.
\end{proof}

\begin{note}\label{rem:conductor}
For $k=1$ the bound \eqref{eq:condbound} always holds: if $n$ is odd then $\Delta\equiv1\pmod4$ and $c$ is odd, so $\Delta_0\equiv1\pmod4$, while if $n$ is even then $4\mid\Delta$, so $c\geq2$. It becomes restrictive for $k\geq2$, as it asks $\Delta$ to have a square factor of size at least $k^2$. It is sharp: for $G=C_3$ and $k=2$ one has $\Delta=48=4^2\cdot3$ and $\mathfrak{f}=12$, so equality holds in \eqref{eq:condbound}, and the fusion ring $C_3+6$ is categorified by the two complex conjugate $C^\ast$ near-group categories of \cite[Theorem 10.19]{MR3635673}.
\end{note}

\subsection{Abelian groups of square order}\label{subsec:2groups}

The levels $\ell$ for which $R(G,\ell)$ is categorifiable are known for $G=C_1$ \cite{MR1981895}, $C_2$ \cite{MR3427429}, $C_3$ \cite[Section 4]{MR3229513} and the elementary abelian $2$-groups \cite[Theorem 1.1]{MR4655273}. We add infinitely many abelian groups of square order (Theorem \ref{thm:2groups}), using the twist identities of Lemma \ref{lem:moments}.

The values $E_r$ often lie in a small field, while the twists are roots of unity, so we may average over the Galois group fixing that field. For a positive integer $L$ let $\zeta_L:=e^{2\pi i/L}$. For a prime $p\equiv1\pmod4$ we put $\sqrt p:=\sum_{t=1}^{p-1}\Jac tp\zeta_p^t$, a square root of $p$ in $\mathbb{Q}(\zeta_p)$.

\begin{lemma}\label{lem:average}
Let $Q\geq2$ be a power of $2$, let $P\geq1$ be an odd squarefree integer whose prime factors are all $\equiv1\pmod4$, and let $L$ be divisible by $Q$ and $P$. Write $\theta\in\mu_L$ as $\theta=\theta_2\,\theta'\prod_{p\mid P}\theta_p$, where $\theta_2$ has $2$-power order, $\theta_p$ has $p$-power order, and $\theta'$ has order prime to $2P$.
\begin{enumerate}
\item Let $\mathbb{K}:=\mathbb{Q}(\zeta_Q,\sqrt p:p\mid P)$, and let $\mathrm{Av}$ be the average over $\mathrm{Gal}(\mathbb{Q}(\zeta_L)/\mathbb{K})$, that is, $\mathrm{Av}=[\mathbb{Q}(\zeta_L):\mathbb{K}]^{-1}\,\mathrm{Tr}_{\mathbb{Q}(\zeta_L)/\mathbb{K}}$. Then
\[
\mathrm{Av}(\theta)=\lambda_\theta\,A(\theta_2)\prod_{p\mid P}B_p(\theta_p).
\]
Here $\lambda_\theta$ is the average of $(\theta')^t$ over $t\in(\mathbb{Z}/\operatorname{ord}\theta')^\times$, i.e.\ the M\"obius function of $\operatorname{ord}\theta'$ divided by its Euler totient; $A(\theta_2)$ is $\theta_2$ if $\theta_2^Q=1$ and $0$ otherwise; and
\[
B_p(1)=1,\qquad B_p(\zeta_p^j)=\frac{-1+\Jac jp\sqrt p}{p-1}\ \ (p\nmid j),\qquad B_p(\theta_p)=0\ \text{ if }\theta_p^p\neq1.
\]
\item Suppose $P>1$, and let $\mathrm{Av}_P$ be the average over $\mathrm{Gal}(\mathbb{Q}(\zeta_L)/\mathbb{Q}(\sqrt P))$. Then the $\sqrt P$-coordinate of $\mathrm{Av}_P(\theta)\in\mathbb{Q}\oplus\mathbb{Q}\sqrt P$ has absolute value at most $1/\varphi(P)$, and it is $0$ if the order of $\theta$ is prime to $P$.
\end{enumerate}
\end{lemma}

\begin{proof}
Write $L=L_2L'\prod_{p\mid P}L_p$, where $L_2$ is a power of $2$, $L_p$ is a power of $p$, and $L'$ is prime to $2P$, so that $\mathrm{Gal}(\mathbb{Q}(\zeta_L)/\mathbb{Q})=(\mathbb{Z}/L_2)^\times\times(\mathbb{Z}/L')^\times\times\prod_{p\mid P}(\mathbb{Z}/L_p)^\times$. The subgroup $\{t\equiv1\bmod Q\}\times(\mathbb{Z}/L')^\times\times\prod_{p\mid P}\{\text{squares}\}$ fixes $\mathbb{K}$, and its index is $\varphi(Q)\,2^{\#\{p\mid P\}}=[\mathbb{K}:\mathbb{Q}]$, since $\mathbb{Q}(\zeta_Q)$ and $\mathbb{Q}(\sqrt p:p\mid P)\subseteq\mathbb{Q}(\zeta_P)$ are linearly disjoint. So this subgroup is $\mathrm{Gal}(\mathbb{Q}(\zeta_L)/\mathbb{K})$, and $\mathrm{Av}(\theta)$ is the product of the averages of $\theta'$, $\theta_2$ and the $\theta_p$ over the corresponding factors. The average of $\theta'$ is $\lambda_\theta$, by definition. The $t\equiv1\bmod Q$ in $(\mathbb{Z}/L_2)^\times$ are the $1+Qs$ with $s\in\mathbb{Z}/(L_2/Q)$, so the average of $\theta_2^t$ is $\theta_2$ times the average of $\theta_2^{Qs}$, which is $\theta_2$ if $\theta_2^Q=1$ and $0$ otherwise, that is, $A(\theta_2)$. Finally, the squares in $(\mathbb{Z}/L_p)^\times$ contain the kernel of reduction modulo $p$, so the average of $\theta_p$ vanishes if $\theta_p^p\neq1$, while for $\theta_p=\zeta_p^j$ with $p\nmid j$ it is $\frac{2}{p-1}\sum_{t\in(\mathbb{F}_p^\times)^2}\zeta_p^{jt}=B_p(\zeta_p^j)$, by the definition of $\sqrt p$. This proves (1). For (2), take $Q=2$ in (1). Then $\mathrm{Av}_P=\mathrm{Av}''\circ\mathrm{Av}$, where $\mathrm{Av}''$ is the average over $\mathrm{Gal}(\mathbb{K}/\mathbb{Q}(\sqrt P))$. Expanding the product in (1), the coefficient of $\sqrt P=\prod_{p\mid P}\sqrt p$ in $\mathrm{Av}(\theta)$ is $\lambda_\theta A(\theta_2)\prod_{p\mid P}\bigl(\pm1/(p-1)\bigr)$ if every $\theta_p$ has order $p$, and $0$ otherwise; in particular it is $0$ if every $\theta_p$ is $1$. The average $\mathrm{Av}''$ fixes $1$ and $\sqrt P$, and kills $\prod_{p\in S}\sqrt p$ for every other nonempty set $S$ of primes dividing $P$. Hence the $\sqrt P$-coordinate of $\mathrm{Av}_P(\theta)$ is that coefficient. As $|\lambda_\theta A(\theta_2)|\leq1$, its absolute value is at most $\prod_{p\mid P}1/(p-1)=1/\varphi(P)$.
\end{proof}

The bound in part (2) is Larson's \cite[Theorem 6.2]{MR3229513}: writing $r+s\sqrt P$ with $r,s\in\mathbb{Q}$ as a sum of roots of unity takes at least $|s|\,\varphi(2P)=|s|\,\varphi(P)$ of them, as $P$ is odd. The new content of (2) is the vanishing for $\theta$ of order prime to $P$, which is what removes the terms $q(g)^2$ in the proof of Proposition \ref{prop:2groups}. Part (1) is used only in the proof of (2).

In Lemma \ref{lem:moments}, the $\sqrt\Delta$-coefficient of $E_2$ is $\nu_2(x)-\frac k2|G[2]|$, which grows with $|G[2]|$, while $E_0\leq2n+\frac12k^2(n^2-1)$. So when $|G[2]|$ is large $E_2$ alone is often enough. We make this precise for abelian groups of square order, where the arithmetic of $\Delta_0$ is rigid. Write $2G:=\{h^2:h\in G\}$, so that $|G[2]|\,|2G|=n$. If $n=4^jw^2$ with $w$ odd, then the odd part of $G$ has order $w^2$ and lies in $2G$, so $|2G|=w^2|2G_2|$, where $G_2$ is the Sylow $2$-subgroup of $G$. The odd part enters the bound of Proposition \ref{prop:2groups} through $n$ and $|2G|$, and it restricts $\Delta_0$ through congruences modulo $w$ (Lemma \ref{lem:2adic}). The elementary abelian $2$-groups, i.e.\ $|2G|=1$, were treated in this way in \cite[Section 4.2]{MR4655273}.

The key input is Lemma \ref{lem:average}(2), applied to $E_2$.

\begin{proposition}\label{prop:2groups}
Let $G$ be an abelian group of order $n=4^jw^2$ with $w$ odd and $j\geq2$, and suppose that $G+kn$ is categorifiable. Put $\tilde k:=2^{j-1}wk$, so that $\Delta=4n(\tilde k^2+1)$, and write $\tilde k^2+1=y^2\Delta_0$. Then $\Delta_0$ is a product of primes $\equiv1\pmod4$ that do not divide $n$, $\Delta_0\equiv1\pmod8$ if $j\geq3$, $\Delta_0$ is a square modulo every prime dividing $w$, and
\begin{equation}\label{eq:2groups}
\varphi(\Delta_0)\,\Bigl|\nu_2(x)-\tfrac k2|G[2]|\Bigr|\,2^{j+1}wy\ \leq\ 2n+\tfrac12k^2(n^2-1).
\end{equation}
In particular, if $n>2|2G|$, then
\[
\varrho(\Delta_0):=\frac{\varphi(\Delta_0)}{\sqrt{\Delta_0}}\ \leq\ \frac{|2G|\,(n+4)}{n-2|2G|}.
\]
\end{proposition}

\begin{proof}
Since $\tilde k$ is even and divisible by every prime dividing $n$, every prime dividing $\tilde k^2+1$ is $\equiv1\pmod4$ and does not divide $n$. So $\Delta_0\equiv1\pmod 4$, $c=2^{j+1}wy$ and $\mathfrak{f}=\Delta_0$. If $j\geq3$, then $4\mid \tilde k$, so $\tilde k^2+1\equiv1\pmod{16}$ and $\Delta_0\equiv y^{-2}\equiv1\pmod8$. If $p\mid w$ is prime, then $p\mid \tilde k$, so $p\nmid y$ and $\Delta_0\equiv y^{-2}\pmod p$.

By Proposition \ref{prop:local}(4), $q(g)=t_g(g)$ is a root of unity whose order divides that of $g$, so it is prime to $\Delta_0$. Average the identity for $E_2$ in Lemma \ref{lem:moments} over $\mathrm{Gal}(\mathbb{Q}(\zeta_N)/\mathbb{Q}(\sqrt{\Delta_0}))$; here $2\Delta_0\mid N$ by Corollary \ref{cor:fsexp}(2). By Lemma \ref{lem:average}(2) with $P=\Delta_0$, the averages of the $q(g)^2$ have $\sqrt{\Delta_0}$-coordinate $0$, and those of the $\theta_W^2$ have $\sqrt{\Delta_0}$-coordinate at most $1/\varphi(\Delta_0)$ in absolute value. So $\bigl|\nu_2(x)-\frac k2|G[2]|\bigr|\,c\,\varphi(\Delta_0)\leq E_0$, and $E_0\leq2n+\frac12k^2(n^2-1)$ by Lemma \ref{lem:moments} and $u_g\leq k^2/4$. This is \eqref{eq:2groups}.

For the last claim, put $a:=|2G|$, so that $|G[2]|=n/a$, and $\varrho:=\varrho(\Delta_0)$. Since $\nu_2(x)=\pm1$ \cite[Theorem 2.7]{MR3427429}, we have $|\nu_2(x)-\frac{kn}{2a}|\geq\frac{kn}{2a}-1>0$. Together with $c>kn/\sqrt{\Delta_0}$, \eqref{eq:2groups} gives $h(k)<0$, where
\[
h(k):=\Bigl(\frac{\varrho}{2a}-\frac12\Bigr)k^2n^2-\varrho kn-2n .
\]
Suppose that $\varrho>a(n+4)/(n-2a)$. Then $\varrho>a$ and $h(1)>0$. Moreover $h(k)-kh(1)=(k^2-k)\bigl(\frac\varrho{2a}-\frac12\bigr)n^2+2(k-1)n\geq0$, so $h(k)>0$ for every $k\geq1$, a contradiction.
\end{proof}

Which $\Delta_0$ occur is governed by the $2$-adic valuation of a Pell solution, and by a congruence modulo $w$.

\begin{lemma}\label{lem:2adic}
Let $j\geq2$, let $w\geq1$ be odd, and let $\Delta_0>1$ be squarefree. Then $\Delta_0$ is the squarefree part of $(2^{j-1}wk)^2+1$ for some $k\geq1$ if and only if the equation $X^2-\Delta_0Y^2=-1$ is solvable in integers, its fundamental solution $(X_0,Y_0)$ has $v_2(X_0)\geq j-1$, and $w$ divides $X$ for some odd power $X+Y\sqrt{\Delta_0}$ of $\varepsilon:=X_0+Y_0\sqrt{\Delta_0}$.
\end{lemma}

\begin{proof}
Suppose $\tilde k^2+1=y^2\Delta_0$ with $\tilde k=2^{j-1}wk$. Then $\tilde k+y\sqrt{\Delta_0}$ is a positive unit of norm $-1$ in $\mathbb{Z}[\sqrt{\Delta_0}]$, so it equals $\varepsilon^{2i+1}$ for some $i\geq0$. Write $\varepsilon^2=U+V\sqrt{\Delta_0}$, so that $U=2X_0^2+1$ is odd and $V=2X_0Y_0$. Also write $\varepsilon^{2i}=U_i+V_i\sqrt{\Delta_0}$. By induction on $i$, $U_i$ is odd and $2^{v_2(X_0)+1}\mid V_i$. Hence $\tilde k=X_0U_i+\Delta_0Y_0V_i$ has $v_2(\tilde k)=v_2(X_0)$. Conversely, if $\varepsilon^{2i+1}=\tilde k+y\sqrt{\Delta_0}$ with $w\mid \tilde k$, then $v_2(\tilde k)=v_2(X_0)\geq j-1$, so $k:=\tilde k/(2^{j-1}w)$ is an integer.
\end{proof}

The coordinates $X$ of the odd powers $X+Y\sqrt{\Delta_0}$ of $\varepsilon$ are periodic modulo $w$, so the last condition is a finite check. It implies that $\Delta_0$ is a square modulo every prime dividing $w$, and it fails whenever a prime $p\mid w$ divides $Y_0$, since then $X\equiv X_0^{2i+1}\not\equiv0\pmod p$ for $\varepsilon^{2i+1}=X+Y\sqrt{\Delta_0}$.

\begin{theorem}\label{thm:2groups}
For every $a\in\mathbb{Z}_{\geq1}$ there is $j_0(a)$ with the following property. Let $G$ be an abelian group of order $4^jw^2$ with $w$ odd, $|2G|=a$ and $j\geq j_0(a)$. Then $G+k|G|$ is not categorifiable for any $k\geq1$, and $R(G,\ell)$ is categorifiable if and only if $\ell=0$. For example, one can take $j_0(2)=3$, that is, the conclusion holds for $G\cong C_4\times C_2^{2j-2}$ with $j\geq3$.
\end{theorem}

\begin{proof}
Put $n=4^jw^2$, and let $T:=a(n+4)/(n-2a)$. Note that $T$ decreases in $n$ and that $n\geq4^j$. By Proposition \ref{prop:2groups}, it suffices to show that $\varrho(\Delta_0)>T$ for every $\Delta_0$ allowed by Proposition \ref{prop:2groups} and Lemma \ref{lem:2adic}. Recall that $\varrho$ is multiplicative, and that $\varrho(p)=(p-1)/\sqrt p>1$ increases with $p$ for odd primes $p$. As the allowed $\Delta_0$ are products of primes $\equiv1\pmod4$, only finitely many of them have $\varrho(\Delta_0)\leq2a$, and each of them has a fixed $v_2(X_0)$ (or no solution of $X^2-\Delta_0Y^2=-1$). By Lemma \ref{lem:2adic}, all of them are excluded once $j-1$ exceeds the largest of these valuations. Moreover $T\leq a(4^j+4)/(4^j-2a)<2a$ once $4^j>4a+4$. This gives $j_0(a)$.

For $a=2$, recall that $|2G|=w^2|2G_2|$, where $G_2$ is the Sylow $2$-subgroup of $G$; so $w=1$ and $G\cong C_4\times C_2^{2j-2}$. Let $j\geq3$. Then $\Delta_0\equiv1\pmod8$ is a product of primes $\equiv1\pmod4$ by Proposition \ref{prop:2groups}, so $\Delta_0\geq17$. As $\varrho(p)$ increases with $p$, a prime $\Delta_0$ has $\varrho(\Delta_0)\geq\varrho(17)=16/\sqrt{17}\approx3.88$, and a composite $\Delta_0$ has $\varrho(\Delta_0)\geq\varrho(5)^2=16/5$. Hence $\varrho(\Delta_0)\geq16/5>34/15\geq T$.

It remains to prove the second assertion: $R(G,\ell)$ is categorifiable if and only if $\ell=0$. The ring $R(G,0)$ is categorified by Tambara--Yamagami categories \cite{TY98}. Now let $\ell>0$. In both cases $n\geq4^j>4a+4>2a$, while a cyclic group has $|2G|\geq n/2$, so $G$ is not cyclic. Hence, by Lemma \ref{lem:levels}, if $R(G,\ell)$ were categorifiable, then $\ell=kn$ for some $k\geq1$, and $R(G,\ell)$ would be the ring $G+kn$, which is not categorifiable by the first assertion.
\end{proof}

\begin{note}\label{rem:2groups}
For every $G$, the identity for $E_2$ in Lemma \ref{lem:moments} can be written as
\[
\sum_{W\in\mathcal{W}}m_W^2\,\theta_W^2+\sum_{g\in G}\Bigl(\frac{k^2}2-2u_g\Bigr)q(g)^2=\Bigl(\nu_2(x)-\frac k2|G[2]|\Bigr)c\sqrt{\Delta_0},
\]
a combination of roots of unity with coefficients in $\frac12\mathbb{Z}_{\geq0}$ of total weight $2n+\frac12k^2n^2$. Larson's bound for arbitrary squarefree $\Delta_0$ \cite[Theorem 6.2]{MR3229513}, applied to twice this identity, gives
\[
\varphi(2\Delta_0)\,\Bigl|\nu_2(x)-\frac k2|G[2]|\Bigr|\,c\ \leq\ 2n+\frac12k^2n^2 .
\]
In the setting of Proposition \ref{prop:2groups}, this is \eqref{eq:2groups} with $\frac12k^2$ more on the right; the proof there removes the $q(g)^2$ by averaging. Along a Pell family of $k$ with fixed $\Delta_0$, the ratio of the two sides tends to $\tilde\varrho(\Delta_0)/|2G|$, where $\tilde\varrho(\Delta_0):=\varphi(2\Delta_0)/\sqrt{\Delta_0}$. So the inequality excludes all but finitely many members of the family if $\tilde\varrho(\Delta_0)>|2G|$, and only finitely many if $\tilde\varrho(\Delta_0)<|2G|$.
\end{note}

\paragraph*{Acknowledgments}
Y. Wang  is supported by NSFC Grant No.\ 12301045, Y. Wang and Z. Yu are supported by NSFC Grant  No.\ 12571041.

The authors acknowledge the use of Claude (Anthropic) in the development of this manuscript. This assistance included suggestions for intermediate claims, proof steps, and corrections, drafting and revision of expository text, literature and citation searches, and computational checks of the examples. All mathematical statements, proofs, computational results, and references were independently checked by the authors, who take sole responsibility for the contents of the paper.

\begin{flushleft}\small
Andrew Schopieray\\
Mathematics Department, SUNY Plattsburgh, Plattsburgh, NY 12901, USA\\
\texttt{ascho013@plattsburgh.edu}

\medskip
Yilong Wang\\
Beijing Institute of Mathematical Sciences and Applications (BIMSA), Beijing, China\\
\texttt{wyl@bimsa.cn}

\medskip
Zhiqiang Yu\\
School of Mathematics, Yangzhou University, Yangzhou 225002, China\\
\texttt{zhiqyumath@yzu.edu.cn}
\end{flushleft}

\end{document}